\documentclass[11pt,a4paper]{amsart}
\usepackage[utf8]{inputenc}
\usepackage{amsmath}
\usepackage{amsfonts}
\usepackage{amssymb}
\usepackage{amsthm}
\usepackage{setspace}
\usepackage{hyperref}
\usepackage{graphicx}
\usepackage{stmaryrd}
\usepackage[shortlabels]{enumitem}
\usepackage{mathtools}
\usepackage{subcaption}
\usepackage[T1]{fontenc}
\usepackage{xcolor}
\usepackage{float}
\usepackage{tikz-cd}
\usepackage{tikz}
\usepackage[noadjust]{cite}
\usepackage{cite}
\usepackage{ragged2e}
\usepackage[margin=1.2in]{geometry}
\usepackage{cleveref}
\newtheorem{thmapp}{Theorem}
\newtheorem{theorem}{Theorem}[section]
\newtheorem{proposition}[theorem]{Proposition}
\newtheorem{lemma}[theorem]{Lemma}
\newtheorem*{mtheorem}{Main Theorem}

\newtheorem{MCorollary}{Corollary}

\usepackage{hyperref}

\usepackage[scr=boondox,scrscaled=1.05]{mathalfa}

\newtheorem{corollary}[theorem]{Corollary}

\newtheorem{case}[theorem]{Case}

\theoremstyle{definition}

\newtheorem{definition}[theorem]{Definition}
\newtheorem{question}[theorem]{Question}

\theoremstyle{remark}
\newtheorem{remark}[theorem]{Remark}

 \newcommand{\Ss}{\mathcal{S}}

 \newcommand{\Th}{\mathrm{Th}}

\newcommand{\ANdim}{\mathop{\rm dim}_\mathrm{AN} }

\newcommand{\U}{\mathcal U}
\newcommand{\V}{\mathcal V}
\newcommand{\cona}{C_1}
\newcommand{\conb}{C_2}
\newcommand{\conc}{C_3}
\newcommand{\cond}{C_0}
\newcommand{\cone}{C_0'}

\theoremstyle{remark}
\newtheorem*{rem-app}{Remark}

 \newcommand{\mr}{\mathrm}

\DeclareMathOperator{\conv}{conv}
\DeclareMathOperator{\diam}{diam}
\DeclareMathOperator{\asdim}{asdim}

\graphicspath{{figures/}}

\title{
Assouad--Nagata dimension of systolic complexes with almost flat intervals\\
}

\newcommand{\anasdim}{\ensuremath{\mathrm{asdim}_\mathrm{AN}}}

\author{Martín Blufstein}
\address[M.~Blufstein]{Departamento de Matemática-IMAS (CONICET), FCEyN, Universidad de Buenos Aires, Buenos Aires, Argentina}
\email{mblufstein@dm.uba.ar}

\author{Victor Chepoi}
\address[V.~Chepoi]{Aix-Marseille Universit\'e and CNRS, LIS, Marseille, France}
\address{Institut Universitaire de France (IUF)}
\email{victor.chepoi@lis-lab.fr}

\author{Huaitao Gui}
\address[H.~Gui]{Institut for Matematiske Fag, University of Copenhagen, 2100 Copenhagen, Denmark}
\email{hg@math.ku.dk}

\author{Damian Osajda}
\address[D.~Osajda]{Institut for Matematiske Fag, University of Copenhagen, 2100 Copenhagen, Denmark}
\address{Instytut Matematyczny
Polskiej Akademii Nauk, ul.\ \'Sniadeckich 8, 00-656 Warszawa, Poland}
\email{dosaj@math.uni.wroc.pl}

\begin{document}

\begin{abstract} 
        We introduce the notion of systolic complexes with almost flat intervals. They constitute a subclass of fit systolic complexes and arise naturally in the study of graphical $C(3)$--$T(6)$ small cancellation complexes.
    	We prove that their Assouad--Nagata dimension is bounded above by a universal constant. As consequences, we obtain universal bounds on the Assouad--Nagata dimension of tame simply connected graphical $C(3)$--$T(6)$ small cancellation complexes and of finitely generated groups defined by graphical $C(3)$--$T(6)$ small cancellation presentations.
        
In the appendix, we observe that the proof of finiteness of asymptotic dimension for buildings due to Dymara and Schick also yields finiteness of their asymptotic Assouad--Nagata dimension.
\end{abstract}

\maketitle

\section{Introduction} 
\label{s:intro}
The \emph{asymptotic dimension} is a coarse invariant of metric spaces introduced by Gromov \cite{Gr}. It plays an important role---among other places---in geometry, group theory, graph theory, and computer science---see surveys \cite{BeDr,Ro}.
The basic question is whether the asymptotic dimension of a given metric space is finite.
A positive answer was given in many important cases (e.g., for Gromov hyperbolic groups \cite{BeDr} and finite dimensional CAT(0) cube complexes \cite{Wr}), but it is still unknown, for example, whether all CAT(0) groups or all finitely presented small cancellation groups have finite asymptotic dimension (see open questions in Section~\ref{s:final}).

The \emph{asymptotic Assouad--Nagata dimension} is an important refinement of the notion of the asymptotic dimension, where additional---linear---control on the diameter of elements of covers is required \cite{DranSmith}.
For uniformly discrete metric spaces the asymptotic Assouad--Nagata dimension coincides with the \emph{Assouad--Nagata dimension} \cite{Ass1} (also known as the \emph{Nagata dimension} \cite{LaSchli2005,GaVeZhe}), and in this article we focus only on such spaces, in particular, on graphs.
Lang and Schlichenmaier \cite{LaSchli2005} gave an important characterization: a metric space has finite Assouad--Nagata dimension if and only if it admits a
quasisymmetric embedding into the product of finitely many trees. See \cite{GaVeZhe} for a very recent application of the finiteness of the Assouad--Nagata dimension for group actions on $\ell_1$.
Within Graph Theory, generalizing and improving a result by Fujiwara and Papasoglu  \cite{FuPa} for planar graphs,  Bonamy, Bousquet, Esperet, Groenland, Liu, Pirot, and Scott  \cite{BoBouEsGrLiPiSc} and Distel \cite{Distel2026} proved that the asymptotic and Assouad--Nagata dimensions of minor-closed graph classes and  of complete Riemannian surfaces are equal to $2$ (see \cite{Liu2025} for a generalization to all minor-closed metrics). 
Assouad--Nagata dimension has also been studied in Theoretical Computer Science, under the name of \emph{weak sparse cover schemes} \cite{AwPe}, or \emph{weak sparse partitions schemes} \cite{JiLiNo+}, and the corresponding graph decompositions are crucial components of several network algorithms.

The asymptotic Assouad--Nagata dimension is an upper bound for the asymptotic dimension, and there are examples of (finitely generated) groups where the difference of the two is arbitrarily large \cite{Nowak2007,Higes2010, BroDyLa2014,Sledd2023}.

\emph{Small cancellation complexes} are $2$--dimensional complexes with some combinatorial nonpositive curvature conditions.
{Small cancellation} is a powerful classical technique for constructing infinite groups.
The \emph{graphical small cancellation} is a far reaching generalization of the classical setting (that is, the one as presented in \cite{LynSch2001,Wise2003}), that has recently allowed for the construction of examples of groups with exotic properties, e.g.\ new nonexact groups \cite{Osa-scl,Os-rf}, or (very recently) groups not coarsely embeddable into every countable family of graphs \cite{Tuck2026}.
There are three usual small cancellation conditions: the $C(6)$ condition, the $C(4)$--$T(4)$ condition, and the $C(3)$--$T(6)$ condition (see Definition~\ref{def:smallcanc}).

The question of finiteness of asymptotic dimension of small cancellation complexes and groups has been open in general; see open questions in Section~\ref{s:final}.
In this article we answer it in the $C(3)$--$T(6)$ case, as a consequence of the following, more general result:

\begin{mtheorem}\label{thm:main}
	Let $X$ be a systolic complex with almost flat intervals. Then the Assouad--Nagata dimension of $X$ is at most $676$.
\end{mtheorem}

\begin{MCorollary}\label{mcor:gc3t6}
The Assouad--Nagata dimension of a tame simply connected graphical $C(3)$--$T(6)$ small cancellation complex is at most $676$.
\end{MCorollary}

The metric on the $C(3)$--$T(6)$ complex considered in the previous corollary is the path metric on its $1$--skeleton (which coincides with a word metric for corresponding groups as in the next corollary). There are other ways of metrizing such complexes---see discussion in Section~\ref{s:c3t6}---and ``tame'' means that they are all quasi-isometric. Assuming tameness is essential; see Remark~\ref{rem:nottame}.

\begin{MCorollary}\label{mcor:c3t6group}
The Assouad--Nagata dimension of a  group defined by a graphical $C(3)$--$T(6)$ small cancellation presentation with finitely many generators is at most $676$.
\end{MCorollary}

In fact, we believe that the Assouad--Nagata dimension of tame simply connected graphical $C(3)$--$T(6)$ small cancellation complexes and finitely generated groups is $338$ times smaller than the bound given above (see Section~\ref{s:final}).
Our methods can be refined to obtain a lower bound, but nowhere near the conjectured bound of $2$.

\emph{Systolic complexes} are flag simplicial complexes whose $1$-skeleta are \emph{bridged} graphs, that is, graphs in which isometric cycles have length $3$.
Bridged graphs play an important role in Metric Graph Theory, and were first introduced by Soltan and Chepoi \cite{SoCh}, and Farber and Jamison \cite{FaJa}.
Systolic complexes were first studied in \cite{Ch_CAT}, and they were rediscovered by Januszkiewicz and \'Swi{\c a}tkowski \cite{JaSw}, and Haglund \cite{Hag2003} within the framework of Geometric Group Theory.
A group is \emph{systolic} if it acts \emph{geometrically}, that is, properly and cocompactly, on a systolic complex.
Systolic complexes and groups should be thought of as nonpositively curved objects and have been successfully used for constructing new exotic examples of groups \cite{JS2003,JaSw}, and for equipping ``classical'' groups with a nonpositive curvature structure, thus allowing to prove new results about them \cite{HoOs, OsPry2018, FukOg2020}. Systolic complexes with \emph{almost flat intervals} is a subclass of \emph{fit systolic complexes} introduced and explored recently in \cite{BlChGuOs}.
With every simply connected  graphical $C(3)$--$T(6)$ small cancellation complex there is associated a systolic complex with almost flat intervals (see Section~\ref{s:almfla}).

The \hyperref[thm:main]{Main Theorem} implies immediately that $2$--dimensional systolic complexes have finite Assouad--Nagata dimension.
Even this was not known before.
In fact, in that case we can get a better bound in Corollary~\ref{cor:2dimsys} (and Theorem~\ref{thm:flat}), with a technically simpler proof in Section~\ref{sec:flat}. 
$2$-dimensional systolic complexes as well as simply connected  graphical $C(3)$--$T(6)$ small cancellation complexes can be equipped with a CAT(0) metric \cite{Gui2025}.
Note that in general $2$--dimensional CAT(0) complexes may  have infinite asymptotic dimension.
Based on a construction by Kondo \cite{Kon2012}, Fukaya and Matsuka \cite{FukMat2026} construct proper $2$-dimensional CAT(0) complexes not coarsely embeddable into a Hilbert space, and hence of infinite asymptotic dimension. 

Because simply connected classical $C(3)$--$T(6)$ small cancellation complexes can be ``systolized'' in a two--dimensional way, we obtain the following result (its proof can be found at the beginning of Section~\ref{sec:flat}). Note that here we do not need to assume tameness of the complex, nor finite generation of the group.

\begin{MCorollary}\label{mcor:clc3t6}
The Assouad--Nagata dimension of a simply connected classical $C(3)$--$T(6)$ small cancellation complex is at most $338$. In particular, the Assouad--Nagata dimension of a  group defined by a classical $C(3)$--$T(6)$ small cancellation presentation (considered with the path metric of the Cayley graph of the presentation) is at most $338$.
\end{MCorollary}

Our proof of the \hyperref[thm:main]{Main Theorem} relies on the characterization of finiteness of the asymptotic dimension  (in fact, also of the asymptotic dimension growth) by Arzhantseva, Niblo, Wright, and Zhang \cite{ArNiWrZh}.
We observe that their results \cite[Corollary 3.2 \& Corollary 3.3]{ArNiWrZh} can be modified to a characterization of finiteness of the asymptotic Assouad--Nagata dimension---Theorem~\ref{t:charANasdimfin} and Corollary~\ref{c:charANasdimfin} \footnote{This was also observed in \cite[page 3, example (B)]{LaStadU2025}}.

In the \hyperref[s:app]{Appendix} we observe that the proof of finiteness of the asymptotic dimension of buildings by Dymara and Schick \cite{DymSchi2009} gives also the analogous bound for the asymptotic Assouad--Nagata dimension. The result is unrelated to the \hyperref[thm:main]{Main Theorem}, except from the fact that there are buildings that are systolic with almost flat intervals, for example, affine $\widetilde{A}_2$ buildings.

\section{Preliminaries} \label{s:prelim}
\subsection{Asymptotic dimension and (asymptotic) Assouad--Nagata dimension} 
For a family $\mathcal U$ of subsets of a metric space $(X,d)$, we define its \emph{mesh} as $\mathrm{mesh}(\mathcal{U}):= \sup \{ \mathrm{diam}(U): \; U\in \mathcal U\}$, and its \emph{$\ell$-multiplicity}, denoted $m_\ell(\mathcal U)$, as the minimal integer $n$ such that for any $x\in X$, the ball $B_\ell(x)$ intersects at most $n$ elements of $\mathcal{U}$.

\begin{definition}
\label{d:asdim} The \emph{asymptotic dimension} of a metric space $X$ is at most $N$, denoted $\asdim X \leq N$ if for every $\ell\in \mathbb{N}$ there exists a cover $\mathcal U_\ell$ of $X$ with finite mesh and $m_\ell(\mathcal{U}_\ell)\leq N+1$.
\end{definition}

\begin{definition}
\label{d:Assouad-Nagata} The \emph{Assouad--Nagata dimension} of a metric space $X$ is at most $N$, denoted $\ANdim X\leq N$, if there is $C>0$ such that for every $\ell\in \mathbb{N}$ there exists a cover $\mathcal U_\ell$ of $X$ with $\mathrm{mesh}(\mathcal{U}_\ell)\leq C\cdot \ell$ and $m_\ell(\mathcal{U}_\ell)\leq N+1$.    
\end{definition}
    
\begin{definition}
\label{d:as-Assouad-Nagata} The \emph{asymptotic Assouad--Nagata dimension} of a metric space $X$ is at most $N$, denoted $\anasdim X\leq N$, if there are $C>0$ and $R_0>0$ such that for every $\ell\in \mathbb{N}$ with $\ell\geq R_0$ there exists a cover $\mathcal U_\ell$ of $X$ with $\mathrm{mesh}(\mathcal{U}_\ell)\leq C\cdot \ell$ and $m_\ell(\mathcal{U}_\ell)\leq N+1$.        
\end{definition}    

Arzhantseva, Niblo, Wright, and Zhang \cite{ArNiWrZh} provide a characterization of finiteness of the asymptotic dimension  (in fact, also of the asymptotic dimension growth) of spaces. They apply it to coarsely median spaces and CAT(0) cubical complexes.
We observe that \cite[Corollary 3.2 \& 3.3]{ArNiWrZh} characterizing finiteness of the asymptotic dimension can be modified to a characterization of finiteness of the asymptotic Assouad--Nagata dimension.
Furthermore, we simplify slightly their criterion: In our Theorem~\ref{t:charANasdimfin} and Corollary~\ref{c:charANasdimfin} below, the index $k$ runs over $0,1,\ldots, \ell$ instead of $1,\ldots,3\ell$ in \cite[Corollary 3.2 \& 3.3]{ArNiWrZh}; and our formulation of 
Theorem~\ref{t:charANasdimfin}(\ref{i:AN2})\ref{i:AN1iii} is simpler than in \cite[Corollary 3.2(2)iii)]{ArNiWrZh}. 
Besides being technically simpler, our formulation allows us to obtain a better bound in the \hyperref[thm:main]{Main Theorem} and in Theorem~\ref{thm:flat} (because we need to use fewer sets $S(x,k,\ell)$).
The proofs are basically the same as the proofs of \cite[Corollary 3.2 \& 3.3]{ArNiWrZh} with only slight quantitative modifications.
For completeness, we present the proofs below.

A metric space $(X,d)$ is \emph{uniformly discrete} if there is a constant $K>0$ such that $d(x,y)>K$ for all $x\neq y$.
We note that for a uniformly discrete metric space $X$ we have that $\ANdim X = \anasdim X$.
Without loss of generality, we can assume that the metric only takes integer values. 

\begin{theorem}
	\label{t:charANasdimfin}
	Let $(X,d)$ be a uniformly discrete metric space, and let $N\in \mathbb{N}$. Then the following are equivalent.
	\begin{enumerate}
		\item \label{i:AN1} $\ANdim X \leq N$;
		\item \label{i:AN2} $\exists D >0$, $\forall \ell \in \mathbb{N}$, $\forall k=0,1,\ldots, \ell$, $\forall x \in X$, we can assign a nonempty subset $S(x,k,\ell) \subseteq X$, satisfying:
		\begin{enumerate}[i)]
			\item $S(x,k,\ell) \subseteq B_{D\cdot \ell}(x)$ for all $\ell \in \mathbb{N}$ and $k=0,\ldots, \ell$;
			\item $\forall \ell \in \mathbb{N}$, $\forall k,k'$ with $0 \leq k \leq k' \leq \ell$, $\forall x \in X$,
			we have $S(x,k,\ell) \subseteq S(x,k',\ell)$;
			\item \label{i:AN1iii} $\forall x,y \in X$ with $d(x,y) \leq \ell$, we have $S(y,k,\ell) \subset S(x,k+d(x,y),\ell)$ for $k=0,\ldots,\ell-d(x,y)$;
			\item $\forall \ell \in \mathbb{N}$, $\forall k=0,1,\ldots,\ell$,  $\forall x \in X$, we have $|S(x,k,\ell)| \leq N+1$.
		\end{enumerate}
	\end{enumerate}
\end{theorem}

\begin{proof}[Proof of Theorem~\ref{t:charANasdimfin}]
The proof follows basically verbatim the one of \cite[Theorem 3.1]{ArNiWrZh}.
As mentioned, for uniformly discrete metric spaces the Assouad--Nagata dimension and the asymptotic Assouad--Nagata dimension coincide, so we can change (\ref{i:AN1}) for $\anasdim X \leq N$.
~\\[-0.5cm]
\medskip

\noindent
(\ref{i:AN1})$\Rightarrow$(\ref{i:AN2}): 
  There is $C>0$ as in Definition~\ref{d:Assouad-Nagata}.
  
  Let $\ell \in \mathbb{N}$ be given. Then there exists a cover $\mathcal U_\ell= \{U_i:i \in I\}$ of $X$ with $\mathrm{mesh}(\mathcal U_\ell)\leq C\cdot \ell$ and $m_\ell(\mathcal{U}_\ell)\leq N+1$.
  
  For every $i\in I$ choose $x_i \in U_i$.
      For $k=0,1,\ldots, \ell$ and $x\in X$, we define
      \[S(x,k,\ell)=\{x_i: B_k(x) \cap U_i \neq \emptyset\}.\]
      We claim that for $D:=1+C$ the condition (\ref{i:AN2}) holds.
      \begin{enumerate}[i)]
        \item If $B_k(x) \cap U_i \neq \emptyset$, we can assume $y \in B_k(x) \cap U_i$. Now $d(y,x_i) \leq \mathrm{mesh}(\mathcal{U}_\ell)$, so $d(x,x_i) \leq k+\mathrm{mesh}(\mathcal{U}_\ell) \leq \ell+C\cdot \ell$. In other words,
            \[S(x,k,\ell) \subseteq B_{\ell+C\cdot \ell}(x)=B_{D\cdot \ell}(x).\]
        \item It is immediate, by our definition of sets $S(x,k,\ell)$.
        \item $\forall x,y \in X$ with $d(x,y) \leq \ell$, $\forall k=0,\ldots,\ell-d(x,y)$, suppose $x_j \in S(y,k,\ell)$, i.e. $B_k(y)\cap U_j\neq \emptyset$, then $B_{k+d(x,y)}(x) \cap U_j \neq \emptyset$, therefore, $x_j\in  S(x,k+d(x,y),\ell)$. So we have:
              \[S(y,k,\ell) \subset S(x,k+d(x,y),\ell).\]
        \item $\forall \ell \in \mathbb{N}$, $\forall k=0,1,\ldots,\ell$, $\forall x \in X$, we have
              \[| S(x,k,\ell)|=| \{x_i: B_k(x) \cap U_i \neq \emptyset\}| \leq m_{\ell}(\mathcal{U}_{\ell}) \leq N+1.\]
      \end{enumerate}

\medskip

\noindent
(\ref{i:AN2})$\Rightarrow$(\ref{i:AN1}): We show that we can take $C=2D$ in Definition~\ref{d:Assouad-Nagata}. For every $\ell \in \mathbb{N}$, let $H=\bigcup\limits_{x\in X}S(x,0,\ell)$. Also, $\forall h \in H$, we define $A_h=\{y: h \in S(y,0,\ell)\}$.
        We define $\mathcal{U}_\ell=\{A_h: h\in H\}$.
        Since $\forall x\in X$ and $\forall h\in S(x,0,\ell)$, we have that $x\in A_h$, then $\mathcal{U}_\ell$ is a cover of $X$.
        Since $S(x,0,\ell) \subseteq B_{D\cdot \ell}(x)$, we know that $d(h,y) \leq D\cdot \ell$ for all $y\in A_h$, which implies $\mathrm{mesh}(\mathcal{U}_\ell) \leq 2D\cdot \ell$.
        
        Finally, let us analyse $m_\ell(\mathcal{U}_\ell)$.
        For all $x \in X$, consider $h\in H$ with $B_\ell(x) \cap A_h \neq \emptyset$. Take $y\in B_\ell(x) \cap A_h$, i.e. $d(y,x) \leq \ell$ and $h \in S(y,0,\ell)$. Now, by the assumptions in Condition (\ref{i:AN2}), we have
        \[S(y,0,\ell) \subseteq S(x,d(x,y),\ell) \subseteq S(x,\ell,\ell).\]
        So $m_\ell(\mathcal{U}_\ell) \leq | S(x,\ell,\ell)| \leq N+1$.
\end{proof}

\begin{corollary}
\label{c:charANasdimfin}
Given a connected graph $X=(V,E)$ with vertices $V$ and edges $E$, and equipped with the edge-path metric $d$, then the following are equivalent:
\begin{enumerate}
  \item \label{i:ANa1} $\anasdim X \leq N$;
  \item \label{i:ANa2} $\exists D >0$, $\forall \ell \in \mathbb{N}$, $\forall k=0,1,\ldots, \ell$, $\forall x \in X$, we can assign a nonempty subset $S(x,k,\ell) \subseteq X$, satisfying:
        \begin{enumerate}[i)]
          \item $S(x,k,\ell) \subseteq B_{D\cdot \ell}(x)$ for all $k=0,\ldots, \ell$; 
          \item $\forall \ell \in \mathbb{N}$, $\forall k,k'$ with $0 \leq k \leq k' \leq \ell$, $\forall x \in X$,
                we have $S(x,k,\ell) \subseteq S(x,k',\ell)$;
          \item $\forall x,y \in X$ with $d(x,y)=1$ (i.e. with $x$ and $y$ connected by an edge), we have:
                $S(y,k,\ell) \subseteq S(x,k+1,\ell)$ for all $k=0,1,\ldots, \ell-1$;
          \item $\forall \ell \in \mathbb{N}$, $|S(x,\ell,\ell)| \leq N+1$.
        \end{enumerate}
\end{enumerate}
\end{corollary}
\begin{proof}
    The proof follows the proof of \cite[Corollary 3.3]{ArNiWrZh} (using \cite[Corollary 3.2]{ArNiWrZh}) adjusted to the proof of our Theorem~\ref{t:charANasdimfin} above.

\medskip

\noindent
(\ref{i:ANa1})$\Rightarrow$(\ref{i:ANa2}): This follows from the implication (\ref{i:AN1})$\Rightarrow$(\ref{i:AN2}) in Theorem~\ref{t:charANasdimfin}. 

\medskip

\noindent
(\ref{i:ANa2})$\Rightarrow$(\ref{i:ANa1}): We proceed exactly as in the proof of (\ref{i:AN2})$\Rightarrow$(\ref{i:AN1}) in Theorem~\ref{t:charANasdimfin}: As in there, we show that we can take $C=2D$ and $R_0=1$ in Definition~\ref{d:Assouad-Nagata}, and we define $H=\bigcup\limits_{x\in X}S(x,0,\ell)$, $A_h=\{y: h \in S(y,0,\ell)\}$, and  $\mathcal{U}_\ell=\{A_h: h\in H\}$ the same way. We get that $\mathcal{U}_\ell$ is a cover of $X$ with $\mathrm{mesh}(\mathcal{U}_\ell) \leq 2D\cdot \ell$.

    For  bounding $m_\ell(\mathcal{U}_\ell)$, for every $x \in X$, consider any $h\in H$ with $B_{\ell}(x) \cap A_h \neq \emptyset$. Take $y\in B_\ell(x) \cap A_h$, i.e. $d(y,x) \leq \ell$ and $h \in S(y,0,\ell)$. By the definition of the edge-path metric $d$, we know that there exists a sequence of vertices $y=y_0,y_1,\ldots,y_k=x$ such that $y_i \in V$ for all $i=0,1,\ldots,k$, $d(y_i,y_{i+1})=1$ for all $i=0,1,\ldots,k-1$, and $k\leq \ell$. Now by the hypothesis, we know:
\[S(y,0,\ell) \subseteq S(y_1,1,\ell) \subseteq S(y_2,2,\ell) \subseteq \cdots \subseteq S(y_k, k, \ell)=S(x,k,\ell)\subseteq S(x,\ell,\ell).\]
So $\{h \in H: B_\ell(x) \cap A_h \neq \emptyset\} \subseteq S(x,\ell,\ell)$, which implies
\[m_\ell(\mathcal{U}_\ell) \leq  |S(x,\ell,\ell)| \leq N+1.\]
\end{proof}

\subsection{Bridged graphs and systolic complexes} All simplicial complexes $X$ considered in this article are connected abstract simplicial complexes. We denote by $V$ the 0-simplices and  by $E$ the $1$-simplices of $X$ and by $G=(V,E)$ the \emph{1-skeleton} of $X$.
All graphs $G=(V,E)$ occurring in this note are $1$-skeleta of simplicial complexes, therefore they are simplicial, undirected and connected. For two distinct vertices $u,v\in V$ we write $u\sim v$ when there is an edge connecting $u$ with $v$ and $u\nsim v$ otherwise. 
The subgraph of $G$ \emph{induced by} a subset $A\subseteq V(G)$ is the graph
$G[A]=(A,E')$ such that $uv\in E'$ if and only if $uv\in E$. 
A \emph{clique} is a complete subgraph of $G$, i.e., a subset of pairwise adjacent vertices. 
The \emph{distance}
$d(u,v)=d_G(u,v)$ between two vertices $u$ and $v$ of a graph $G$ is the
length of a shortest $(u,v)$-path. 
An induced subgraph $H=G[A]$ of a graph $G$ is an \emph{isometric subgraph} of $G$ if $d_H(u,v)=d_G(u,v)$ for any
two vertices $u,v\in A$. The \emph{interval}
$I(u,v)$ between $u$ and $v$ consists of all vertices on shortest
$(u,v)$-paths, that is, of all vertices (metrically) \emph{between} $u$
and $v$: $I(u,v)=\{ x\in V: d(u,x)+d(x,v)=d(u,v)\}$. For two vertices $u, v\in V$  with $d(u,v)=\ell$ and $0\leq k\leq \ell$ 
the \emph{$k$-level} $L_k(u,v)$ of the interval $I(u,v)$ is the subgraph of $G$ induced by all  $x\in I(u,v)$ such that $d(u,x)=k$ and $d(x,v)=\ell-k$. An induced
subgraph of $G$ is called \emph{convex}
if it contains the interval in $G$ between any two of its
vertices.  For a vertex $v$ of $G$ and an integer $r\geq 1$, we will denote  by $B_r(v)$ the \emph{ball} in $G$
(and the subgraph induced by this ball)  of radius $r$ centred at  $v$, i.e., $B_r(v)=\{ x\in V: d(v,x)\leq r\}.$ More generally, the $r$-{\it ball  around a set} $A\subseteq V$
is the set (or the subgraph induced by) $B_r(A)=\{ v\in V: d(v,A)\leq r\},$ where  $d(v,A)=\mbox{min} \{ d(v,x): x\in A\}$.
Similarly, the \emph{sphere} in $G$ centred at $v$ of radius $r$ is the set (or the subgraph induced by) $S_r(v)=\{ x \in V(G): d(v,x) = r\}$.
As usual, $N(v)=S_1(v)$ denotes the set of neighbours of a vertex
$v$ in $G$, and $N[v]=B_1(v)$.
The \emph{link} of $v \in V$ is the subgraph of $G$ induced by $N(v)$.  
For a set $K$ of vertices of a graph
$G = (V,E)$ and an integer $k \ge 0$, let
$B^*_k(K)=\bigcap_{s\in K} B_k(s)$.  If $K \subseteq K'$, then $B^*_k(K')\subseteq B^*_k(K)$. If $K$ is a clique, then $B^*_1(K)$ is the union of $K$ and the set of vertices adjacent to all vertices in $K$.

\begin{definition}[Bridged/systolic]
    A graph $G$ is \emph{bridged} if any isometric cycle of $G$ has size $3$. A flag simplicial complex $X$ is \emph{systolic} if its $1$-skeleton is a bridged graph.
\end{definition}

Systolic complexes are exactly simply connected flag simplicial complexes in which the links of vertices are \emph{6-large}, i.e., do not contain induced $4$-cycles and $5$-cycles \cite{Ch_CAT, JaSw}. Bridged graphs have been characterized in the following way: 

\begin{theorem}\label{bridged-convex-balls}  \cite{FaJa,SoCh} A graph is bridged if and only if   the balls around convex sets are convex. 
\end{theorem} 

Each level set $L_k(u,v)$ of an interval $I(u,v)$  of a systolic complex is convex  because $L_k(u,v)$ is the intersection of two convex sets $B_k(u)$ and $B_{\ell-k}(v)$ (where $\ell=d(u,v)$). 

A graph $G$ satisfies the \emph{triangle condition} if for  $v,x,y\in V$ with $d(v,x)=d(v,y)=k$ and $x\sim y$, there exists $z\sim x,y$ with $d(v,z)=k-1$. All bridged graphs satisfy the triangle condition \cite{Ch_mt}.

\subsection{Almost flat intervals in systolic complexes}  \label{s:almfla}
In this paper, by $X$ we will denote a systolic complex and by $G=(V,E)$ the bridged graph which is the $1$-skeleton of $X$.
With some abuse of notation, we will denote by $I(u,v)$ the interval between $u$ and $v$ in $G$ (viewed as a set of vertices), the subgraph of $G$ induced by $I(u,v)$, and the subcomplex of $X$ spanned by $I(u,v)$. We will consider the following type of intervals in systolic complexes:

\begin{definition}[Almost flat intervals in systolic complexes] An interval $I(u,v)$ of a systolic complex $X$ is called \emph{almost flat}  if every level $L_k(u,v)$ of $I(u,v)$ is a line. A \emph{systolic complex with almost flat intervals} is a systolic complex in which all intervals are almost flat. 
\end{definition}

The simplest example of a systolic complex with almost flat intervals is the systolic plane (the subdivision of ${\mathbb R}^2$ into equilateral triangles with sides 1).
Each interval $I(u,v)$ of the resulting triangular grid is a subgrid bounded by a parallelogram having $u$ and $v$ as opposite vertices and  with sides of lengths $k_1$ and $k_2$ such that $k_1+k_2=d(u,v)$. 
In the case of $2$-dimensional systolic complexes ($K_4$-free bridged graphs), it was shown in  \cite{ChLaRa} that each interval $I(u,v)$ is isometric to a subgraph of such a parallelogram of the triangular grid, namely it is isometric to the region bounded by two shortest paths of the triangular grid between the images of $u$ and $v$.
In \cite{ChLaRa}, such intervals have been called \emph{burned lozenges}. See Figures~\ref{fig:wedge} \& \ref{fig:net} for examples.
We will say that a systolic complex $X$ has \emph{flat intervals} if each interval $I(u,v)$ of $X$ is an at most  $2$-dimensional systolic complex, i.e., it is a burned lozenge.

\begin{remark}\label{rem:afiwrtv}
    In \hyperref[thm:main]{Main Theorem} and Theorem~\ref{thm:flat} we show that a systolic complex $X$ with almost flat intervals has finite asymptotic dimension (bounded by a global constant independent of $X$).
    However we note that for either statement to be true, it suffices that $X$ contains a vertex $x_0$ such that intervals of the form $I(x_0, \cdot)$ are (almost) flat.
\end{remark}

\begin{lemma}
\label{lem:almost flat to flat}
     If an interval $I(u,v)$ of a systolic complex $X$ contains no 4-cliques, then $I(u,v)$ is flat. 
\end{lemma}

\begin{proof}
Each interval $I(u,v)$ of $X$ is a simply connected and 6-large subcomplex of $X$, thus $I(u,v)$ is systolic.
Since $I(u,v)$ does not contain 4-cliques, $I(u,v)$ is $2$-dimensional and therefore flat.
\end{proof}

Let $X$ be a systolic complex and $G$ be its $1$-skeleton.
For two vertices $u,v$ with $d(u,v)=\ell$ and $0\leq k_1 \leq k_2 \leq \ell$, let $L_{k_1,k_2}(u,v)$ denote the subgraph of $G$ induced by all $x\in I(u,v)$ such that $k_1 \leq d(u,x)\leq k_2$.
Equivalently, $L_{k_1,k_2}(u,v)$ is the subgraph induced by the vertices lying in the levels $L_{k_1}(u,v),\ldots,L_{k_2}(u,v)$.
We will call $L_{k_1,k_2}(u,v)$ a $(k_1,k_2)$-\emph{slice} of $I(u,v)$, or simply a \emph{slice} of $I(u,v)$.
Furthermore, we classify the edges of $I(u,v)$ as follows: an edge $xy$ of $I(u,v)$ is called \emph{horizontal} if $x,y$ belong to the same level $L_i(u,v)$ and \emph{oblique} if $x$ and $y$ belong to two consecutive levels of $I(u,v)$.

Now, we would like to specify the structure of almost flat intervals $I(u,v)$ in systolic complexes. Namely, we distinguish two types of slices in $I(u,v)$ and we show that $I(u,v)$ is obtained by amalgamating such slices.
A \emph{chopped lozenge} is either a burned lozenge $I(u,v)$ or a burned lozenge $I(u,v)$ without one or two of its ends $u$ and $v$.
A \emph{zip} is a sequence of 4-cliques $K_1,\ldots,K_m$ such that any two non-consecutive cliques are disjoint and the cliques $K_{i-1},K_{i+1}$ intersect $K_i$ in two non-incident edges of $K_i$. 
We say that an interval $I(u,v)$ with $d(u,v)=\ell$  is an \emph{amalgamation} of chopped lozenges and zips if there exist 
 $0 = k_0< \ldots <k_m = \ell$ satisfying the following conditions:
 \begin{itemize}
     \item[(1)] the first slice $L_{k_0,k_{1}}(u,v)$ and the last slice $L_{k_{m-1},k_{m}}(u,v)$ are chopped lozenges;
     \item[(2)] each slice $L_{k_i,k_{i+1}}(u,v)$ is either a chopped lozenge or a zip. When moving from $u$ to $v$, chopped lozenges and  zips alternate, i.e., if $L_{k_i,k_{i+1}}(u,v)$ is a chopped lozenge (respectively, zip), then $L_{k_{i+1},k_{i+2}}(u,v)$ is a zip (respectively, chopped lozenge); 
     \item[(3)] a chopped lozenge  (respectively, zip) $L_{k_{i-1},k_i}(u,v)$  and a zip (respectively, chopped lozenge) $L_{k_i,k_{i+1}}(u,v)$ are merged together along the level $L_{k_i}(u,v)$, which is a single horizontal edge. 
 \end{itemize}
See Figure~\ref{fig:wedge-almostflat} (on the right) for an example of an interval being an {amalgamation} of chopped lozenges and zips.

\begin{lemma} \label{l:structure-almost-flat-intervals}
    If an interval $I(u,v)$ of a systolic complex $X$ is almost flat, then $I(u,v)$ is an amalgamation of chopped lozenges and zips.
\end{lemma}

\begin{proof} First, observe that each clique of $I(u,v)$ must belong to  two consecutive levels of $I(u,v)$. Since each level of $I(u,v)$ is a path, $I(u,v)$ cannot contain 5-cliques. Now, pick any 4-clique $K=\{ x_1,x_2,x_3,x_4\}$. Then necessarily $K$ must belong to two consecutive levels  $L_{k}(u,v)$ and $L_{k+1}(u,v)$, two vertices in one level and two vertices in the second level, say $x_1,x_2\in L_k(u,v)$ and $x_3,x_4\in L_{k+1}$. We assert that $L_k(u,v)=\{ x_1,x_2\}$ and $L_{k+1}(u,v)=\{ x_3,x_4\}$. Suppose by way of contradiction that $L_k(u,v)$ contains a third vertex $z\ne x_1,x_2$. Since $L_k(u,v)$ is a path, we can suppose that $z$ is adjacent to exactly one of $x_1,x_2$, say $z\sim x_2, z\nsim x_1$. By applying the triangle condition to $v$ and the edge $x_2z$, there exists a vertex $y\in L_{k+1}(u,v), y\sim x_2,z$. Since $L_k(u,v)$ is convex, $y$ cannot coincide with $x_3,x_4$. Since $x_2\sim x_3,x_4,y$ and $L_{k+1}(u,v)$ is convex, $y\sim x_3,x_4$, contrary to the assumption that $L_{k+1}(u,v)$ is a path. Consequently, $L_{k,k+1}(u,v)$ is the 4-clique $K=\{ x_1,x_2,x_3,x_4\}$.

Consider all pairs of indices $0<k_i<k_j<\ell$ such that the slice $L_{k_i,k_j}(u,v)$ is a maximal by inclusion zip.
Pick two consecutive such zips $L_{k_i,k_j}(u,v)$ and $L_{k_{i'},k_{j'}}(u,v)$ with $i<j<i'<j'$. Let $L_{j-1}(u,v)=\{ x_{j-1},y_{j-1}\}, L_{j}(u,v)=\{ x_j,y_j\}$ and $L_{i'}(u,v)=\{x_{i'},y_{i'}\}, L_{i'+1}(u,v)=\{ x_{i'+1},y_{i'+1}\}$. 
Since $\{ x_{j-1},y_{j-1}, x_j,y_j\}$ and $\{x_{i'},y_{i'}, x_{i'+1},y_{i'+1}\}$ are cliques, from the choice of these vertices we conclude that the slice $L_{j,i'}(u,v)$ coincides with the interval $I(x_{j-1},x_{i'+1})$ minus the vertices $x_{j-1},x_{i'+1}$ (the same conclusion holds for other three intervals $I(x_{j-1},y_{i'+1}), I(y_{j-1},x_{i'+1})$ and $I(y_{j-1},y_{i'+1})$).
By Lemma \ref{lem:almost flat to flat}, the slice $L_{j,i'}(u,v)$ is a chopped lozenge, which is amalgamated with the zips  $L_{k_i,k_j}(u,v)$ and $L_{k_{i'},k_{j'}}(u,v)$ along the horizontal edges $L_j(u,v)=\{ x_j,y_j\}$ and $L_{i'}(u,v)=\{ x_{i'},y_{i'}\}$.
Finally, if the  zip $L_{i',j'}$ does not exist then we apply a similar argument with $v$ instead of $x_{i'+1}$.
Analogously we can also show that the slice preceding a zip $L_{k_i,k_j}(u,v)$ is a chopped lozenge. 
This shows that every almost flat interval $I(u,v)$ is an amalgamation of chopped lozenges and zips.
\end{proof}

\subsection{Graphical $C(3)$--$T(6)$ small cancellation complexes}
\label{s:c3t6}
A map $X\to Y$ between CW complexes is \emph{combinatorial} if its restriction to every open cell of $X$ is a homeomorphism onto an open cell of $Y$.
A CW complex is \emph{combinatorial} if its attaching maps are combinatorial.
An \emph{immersion} is a locally injective combinatorial map.

Let $\Theta$ be a connected graph without vertices of degree $1$, and let $f\colon \Omega \to \Theta$ be an immersion of graphs. 
Decompose $\Omega = \coprod_i \Omega_i$, where each $\Omega_i$ is a connected component.

For these data, we construct a combinatorial $2$-complex as follows.
Start with $\Theta$ as the $1$-skeleton.
For every immersed cycle $C\rightarrow \Omega$, attach a $2$-cell to $\Theta$ along the composition $C\rightarrow \Omega \rightarrow \Theta$.
The resulting complex is called the \emph{graphical complex} associated to $f\colon \Omega \to \Theta$, and we denote it by $X_{(\Omega,\Theta)}$.

Let $\gamma_1: P\to \Omega_i$ and $\gamma_2\colon P\to \Omega_j$ be two lifts of the same path $P\rightarrow \Theta$, where $i$ and $j$ may be the same. 
We say that these lifts are \emph{essentially distinct} if there is no isomorphism $\Omega_i \rightarrow \Omega_j$ such that the following diagram commutes:
\[
\begin{tikzcd}
    P \arrow[r,"\gamma_1"] \arrow [d,swap,"\gamma_2"] & \Omega_i \arrow[d]\arrow[ld,dashed] \\ 
    \Omega_j \arrow[r] & \Theta
\end{tikzcd}
\]
A \emph{piece} in $X_{(\Omega,\Theta)}$ is a nontrivial immersed path $\gamma: P\rightarrow \Theta$ which admits two essentially distinct lifts to $\Omega$.

A \emph{disc diagram} $\varphi:D\to X_{(\Omega,\Theta)}$ is a combinatorial map, where $D$ is obtained from a combinatorial structure of $S^2$ by removing a $2$-cell.
A \emph{piece} in a disc diagram is a nontrivial path $P\rightarrow D$ that factors in two distinct ways through $2$-cells, i.e., there are $2$-cells $R_1\rightarrow D$ and $R_2\rightarrow D$ such that $P\rightarrow D$ factors both as $P\rightarrow R_1 \rightarrow D$ and $P\rightarrow R_2 \rightarrow D$, and there is no isomorphism $R_1\rightarrow R_2$ making the following diagram commute:
\[
\begin{tikzcd}
    P \arrow[r] \arrow[d] & R_1 \arrow[d] \arrow[ld,dashed]\\
    R_2 \arrow[r] & D
\end{tikzcd}
\]

For a piece $P\rightarrow D$ with such factorizations, each composition $P\rightarrow R_i\rightarrow D\xrightarrow{\varphi} X_{(\Omega,\Theta)}$ determines a canonical lift $\gamma_i: P\rightarrow \Omega_{k_i}$ to the component of $\Omega$ corresponding to $R_i$. 
We say that $\varphi$ is \emph{reducible} along the piece $P\rightarrow D$ if the lifts $\gamma_1$ and $\gamma_2$ are not essentially distinct.
A disc diagram is \emph{reduced} if it is not reducible along any piece.

\begin{definition}
    \label{def:smallcanc}
    For $p,q\in \mathbb{N}^+$, we say that the graphical complex $X_{(\Omega,\Theta)}$ associated to $f:\Omega \to \Theta$ satisfies 
    \begin{itemize}
        \item the \emph{C(p) condition} if there is no immersed cycle $C\rightarrow\Omega$ such that $C\rightarrow \Omega \rightarrow \Theta$ is the concatenation of less than $p$ pieces;
        \item the \emph{T(q) condition} if for every reduced disc diagram $D\rightarrow X_{(\Omega,\Theta)}$, there is no interior vertex of $D$ whose degree is greater than $2$ and less than $q$. 
    \end{itemize}  
\end{definition}

Given a simply connected graphical $C(3)$--$T(6)$ complex $X_{(\Omega,\Theta)}$ we construct a systolic complex as follows.
For each $\Omega_i$ add to $\Theta$ every possible edge so that the vertices in $f(\Omega_i)$ form a clique, and denote the new simplicial graph by $G_{X_{(\Omega,\Theta)}}$.
The flag completion of $G_{X_{(\Omega,\Theta)}}$ is called the \emph{thickening of $X_{(\Omega,\Theta)}$} and denoted by $\Th(X_{(\Omega,\Theta)})$.
It is shown in \cite[Propositions 2.12 \& 2.13] {BlChGuOs} that $\Th(X_{(\Omega,\Theta)})$ is a systolic complex with almost flat intervals.
We say that $X_{(\Omega,\Theta)}$ is \emph{tame} if it is quasi-isometric to $\Th(X_{(\Omega,\Theta)})$. A natural example of a tame $X_{(\Omega,\Theta)}$ is a one with finitely many isomorphism classes of the graphs $\Omega_i$.

\begin{remark}\label{rem:conedoff}
    Yet another $2$--complex associated with a graphical $C(3)$--$T(6)$ complex $X_{(\Omega,\Theta)}$ can be obtained by glueing simplicial cones over relators $\Omega_i$ to $\Theta$. It is shown in \cite{Gui2025} that such a complex carries a CAT(0) metric. If $X_{(\Omega,\Theta)}$ is tame then it is quasi-isometric also to this CAT(0) space.
\end{remark}

A \emph{graphical presentation} is the special case of the preceding construction in which $\Theta =R_S$, where $R_S$ is a wedge of oriented loops, one labeled by each element of a set $S$.
We write such a presentation as $\mathcal{P}=\langle S:f\rangle$ with $f:\Omega=\coprod_{i\in I}\Omega_i \rightarrow R_S$, and assume that every $\Omega_i$ is finite. Equivalently, the map $f$ may be viewed as an $S\cup S^{-1}$-labeling of the oriented edges of $\Omega$, with opposite orientations receiving inverse labels. 

The group defined by $\mathcal{P}$ is $\Gamma=\Gamma(\mathcal{P}):=\pi_1(R_S)/\langle\!\langle f_*(\pi_1(\Omega_i)_{i\in I})\rangle\!\rangle$.
In other words, $\Gamma$ is the quotient of the free group $F(S)$ by the normal closure of all words read along closed paths in $\Omega$.
If $X=X(\mathcal{P})$ denotes the graphical complex associated to $f:\Omega \rightarrow R_S$, then $\Gamma = \pi_1(X)$.
We say that $\mathcal{P}$ is a \emph{$C(p)$--$T(q)$ graphical small cancellation presentation} if the associated graphical complex $X$ is $C(p)$--$T(q)$.
The group $\Gamma$ is \emph{finitely generated} if the set $S$ of generators is finite, and \emph{finitely presented} if the graph $\Omega$ is finite and the set $S$ of generators is finite.

Let $\widetilde{X}=\widetilde{X}(\mathcal{P})$ be the universal cover of $X$. 
Its $1$-skeleton is naturally identified with the Cayley graph $\operatorname{Cay}(\Gamma ,S)$. Since the image of $\pi_1(\Omega_i)$ is trivial in $\Gamma$, every restriction $f|_{\Omega_i}: \Omega_i\rightarrow R_S$ admits lifts to $\operatorname{Cay}(\Gamma,S)$. 
Let $\tilde{f}: \widetilde{\Omega} \rightarrow \operatorname{Cay}(\Gamma,S)$ be the disjoint union of all such lifts. 
The graphical complex associated to $\tilde{f}$ is naturally identified with $\widetilde{X}$.
Two lifts $\tilde{f}_1:\Omega_i\rightarrow \operatorname{Cay}(\Gamma,S)$ and $\tilde{f}_2:\Omega_j\rightarrow \operatorname{Cay}(\Gamma,S)$ are called equivalent if there is an isomorphism $\alpha: \Omega_i\rightarrow \Omega_j$ such that $\tilde{f}_1=\tilde{f}_2\circ \alpha$. 
Retaining one representative from each equivalence class gives a graphical complex denoted by $\widetilde{X}^+$. 
If $X$ is $C(p)$--$T(q)$, then so is $\widetilde{X}^+$; see, for example, \cite[Lemma 3.7]{DG2026} for a proof.

\subsection{Classical small cancellation}
\label{s:classc}

Every combinatorial $2$-complex $X$ can be viewed as a subcomplex of a graphical complex. More precisely, let $X^{(1)}$ be the $1$-skeleton of $X$ and $\{f_\alpha: \partial R_\alpha\rightarrow X^{(1)}\}$ be the collection of the attaching maps of all the 2-cells in $X$. Then $X$ is a subcomplex of the graphical complex associated to $\coprod_\alpha f_\alpha: \coprod_\alpha \partial R_\alpha\rightarrow X^{(1)}$. 
We say that $X$ satisfies the \emph{classical $C(p)$ or $T(q)$ condition} if its associated graphical complex satisfies the corresponding condition.
For a standard presentation of classical small cancellation theory, we refer to Lyndon--Schupp \cite[Chapter V]{LynSch2001} and McCammond--Wise \cite{MW2002}.
 
The following lemma records some standard properties of classical $C(3)$--$T(6)$ complexes.
They follow from the classical disc-diagram theory of Lyndon--Schupp (see in particular \cite[Chapter V, \S2 and \S4]{LynSch2001}), together with Pride's observation that every piece in a $T(q)$-complex, $q\geq 5$, consists of a single edge \cite{Pride1998}.  

\begin{lemma}\label{lem:edge-pieces}
    Let $X$ be a simply connected classical $C(3)$--$T(6)$ complex. Then:
    \begin{itemize}
        \item every 2-cell $R$ of $X$ is embedded;
        \item the intersection of two distinct $2$-cells $R_1$ and $R_2$ is either empty, a vertex or an edge.
    \end{itemize}
\end{lemma}

\begin{proposition}\label{p:c3t6systolic}
    Let $X$ be a simply connected classical $C(3)$--$T(6)$ small cancellation complex. Then there exists a $2$-dimensional systolic complex $Y$ that is quasi-isometric to $X$.
\end{proposition}
\begin{proof}
    We construct $Y$ by subdividing the $2$-cells of $X$ as follows.
    Take a $2$-cell $F$ of $X$ whose boundary has $n$ edges.
    If $n=1$, remove $F$ and its boundary $1$-cell.
    If $n=2$, identify the two edges and remove the interior of the $2$-cell.
    If $n=3$, do nothing.
    If $n=4$, subdivide $F$ into two triangles by adding a diagonal.
    If $n=5$, subdivide $F$ into three triangles by adding two diagonals.
    If $n\geq 6$, express $n$ as $6k+r$, with $0\leq r \leq 5$.
    Construct a tilling of a regular hexagon of side length $k$ and subdivide one of its boundary triangles by adding $r$ edges as shown in Figure \ref{fig:subdivision}.
    \begin{figure}
        \centering
        \includegraphics[width=0.6\linewidth]{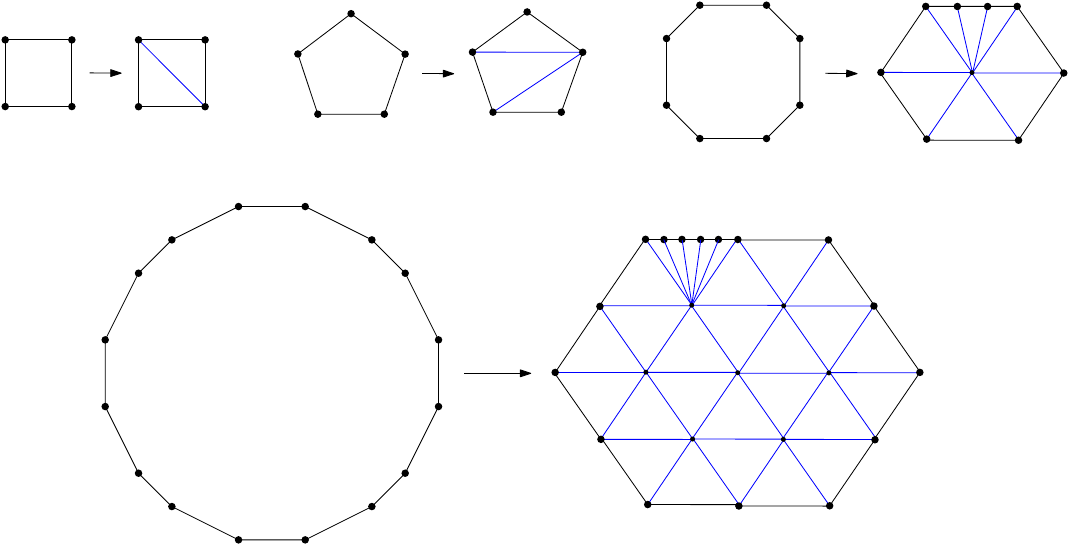}
        \caption{Subdivision of 2-cells with $4$, $5$, or at least $6$ boundary edges resulting in systolic discs}
        \label{fig:subdivision}
    \end{figure}
    Replace the $2$-cell $F$ by this hexagon, and repeat this procedure for all $2$-cells with at least six boundary edges in $X$.
    The resulting complex is $Y$.
    Since $2$-cells embed in $C(3)$--$T(6)$ complexes, $Y$ is a simplicial complex. Since $X$ is simply connected, $Y$ is also simply connected. 
    Using the fact that the intersection between any two $2$-cells in a $C(3)$--$T(6)$ complex is empty, a vertex or a single edge, it is a straightforward check to see that links of vertices in $Y$ are $6$-large, therefore $Y$ is systolic.
    By construction it is quasi-isometric to $X$.
\end{proof}

\section{The flat case}\label{sec:flat}

The goal of this section is to prove the \hyperref[thm:main]{Main Theorem} for systolic complexes with flat intervals. The argument is technically simpler here and will also be used in Section~\ref{s:almflatproof} to prove the general case.
In fact, we obtain a slightly better bound and prove the following.

\begin{theorem}\label{thm:flat}
	Let $X$ be a systolic complex with flat intervals. Then the Assouad--Nagata dimension of $X$ is at most $338$. 
\end{theorem}

As mentioned before, $2$-dimensional systolic complexes, were shown in \cite{ChLaRa} to have flat intervals.
Thus we obtain the following corollary.

\begin{corollary}\label{cor:2dimsys}
    Let $X$ be a $2$-dimensional systolic complex. Then the Assouad--Nagata dimension of $X$ is at most $338$.
\end{corollary}

\begin{proof}[Proof of Corollary~\ref{mcor:clc3t6}]
    It follows from Proposition~\ref{p:c3t6systolic} and Theorem~\ref{thm:flat}. 
\end{proof}

\begin{figure}
    \centering
    \includegraphics[width=0.25\linewidth]{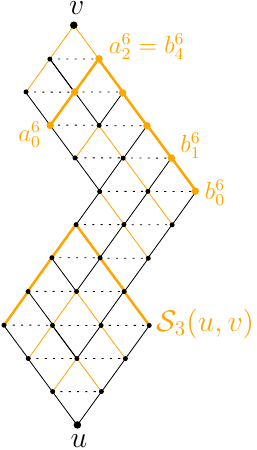}
    \caption{A flat interval.}
    \label{fig:wedge}
\end{figure}

We prove Theorem~\ref{thm:flat} by applying the characterization of Assouad--Nagata dimension provided by Corollary \ref{c:charANasdimfin}.
Let $u,v$ be two vertices of a systolic complex $X$ with flat intervals. 
Then the interval $I(u,v)$ can be embedded in the systolic plane as a burned lozenge. 
The edges belonging to the level sets of $I(u,v)$ are called horizontal.
Now remove all horizontal edges from $I(u,v)$ to obtain a square complex $I^{\diamond}(u,v)$: the vertices of $I^{\diamond}(u,v)$ are the vertices of $I(u,v)$, the edges of $I^{\diamond}(u,v)$ are the oblique edges of $I(u,v)$, and the squares of $I^{\diamond}(u,v)$ correspond to pairs of triangles of $I(u,v)$ sharing a common horizontal edge, see Figure~\ref{fig:wedge}. 
The square complex  $I^{\diamond}(u,v)$ is simply connected and non-positively curved, thus $I^{\diamond}(u,v)$ is CAT(0). 
For any $z\in I(u,v)$, we denote by $d_{nor}(u,z)$ the normal cube path distance between $u$ and $z$ in $I^{\diamond}(u,v)$ \cite{NiRe}.
Let
\[\Ss_k(u,v)=\{ z\in I(u,v): d_{nor}(u,z)=k\}\]
be the sphere of radius $k$ and centre $u$ in $I(u,v)$ with respect to the normal cube path distance.
We call  $\Ss_k(u,v)$ the \emph{$k$-wedge} of $I(u,v)$.
Note that each wedge $\Ss_k(u,v)$ consists of at most two convex paths
\[P_1^k(u,v) = (a_0^k,\ldots,a_{n(k)}^k), \ \ P_2^k(u,v) = (b_0^k,\ldots,b_{m(k)}^k)\]
with a common end-vertex $b_{m(k)}^k=a_{n(k)}^k$, called the \emph{pointer} of  $\Ss_k(u,v)$, see Figure \ref{fig:wedge}. 
The other end-vertices $a_0^k$ and $b_0^k$  are called the \emph{ends} of the wedge  $\Ss_k(u,v)$. 
The \emph{pointer} of $\Ss_k(u,v)$ is the unique vertex of $\Ss_k(u,v)$ closest to $v$ (with respect to the distance of $X$). 
If  $\Ss_k(u,v)$ consists of a single convex path $P_i^k(u,v)$ with $i\in \{ 1,2\}$, we still call the vertex of  $P_i^k(u,v)$ closest to $v$ the \emph{pointer} of $\Ss_k(u,v)$.

The following lemma is central in the proof and is an analogue of \cite[Lemma 5.1]{ArNiWrZh} used in the CAT(0) cubical case. The proof follows the general idea from \cite{ArNiWrZh} using the fact that in our case each interval $I(u,v)$ can be equipped with a structure of a CAT(0) square complex $I^{\diamond}(u,v)$. 

\begin{lemma}\label{l:nets-2d}
Let $X$ be a systolic complex with flat intervals and let $x_0\in X$ be a vertex. Then for all $\ell \in \mathbb{N}$ and $x \in X$, there are subsets $C_x \subseteq I(x_0,x)$, and  maps $p_x\colon I(x_0,x) \rightarrow \mathcal{P}(C_x)\setminus \{ \emptyset\}$, where $\mathcal{P}(C_x)$ denotes the power set of $C_x$, satisfying:
\begin{enumerate}[(a)]
  \item If $d(x,y)=1$ and $y \in I(x_0,x)$, then $C_x \cap I(x_0,y)=C_y$, and $p_x\vert_{I(x_0,y)}=p_y$;
  \item For all $z \in I(x_0,x)$ and $w \in p_x(z)$, we have $d(z,w) \leq 3\ell$; 
  \item For any $z\in I(x_0,x)$, we have $|B_{6\ell}(z) \cap C_x| \leq 338$. 
\end{enumerate}
\end{lemma}

\begin{proof}
    We start by defining the sets $C_x$ for every $x\in X$ and $\ell \in \mathbb{N}$.
    Consider the wedges $\Ss_{0}(x_0,x),\Ss_{\ell}(x_0,x),\Ss_{2\ell}(x_0,x),\ldots,\Ss_{m\ell}(x_0,x)$ where $m=\lfloor\frac{d_{nor}(x_0,x)}{\ell}\rfloor$.

    Recall that an $\ell$-\emph{net} for a set $R$ is a subset $R_0$ of $R$ such that each vertex of $R$ is at distance at most $\ell$ from $R_0$.
    We define an $\ell$-net $N_{i\ell}(x_0,x)$ for each  wedge $\Ss_{i\ell}(x_0,x)$.
    If $\Ss_{i\ell}(x_0,x)$ is a single vertex, then $N_{i\ell}(x_0,x)=\Ss_{i\ell}(x_0,x)$, and if not, then 
    \[N_{i\ell}(x_0,x) = \left\{a_0^{i\ell}, a_\ell^{i\ell},\ldots,a_{\lfloor\frac{n(i\ell)}{\ell}\rfloor\ell}^{i\ell},b_0^{i\ell},b_\ell^{i\ell},\ldots,b_{\lfloor\frac{m(i\ell)}{\ell}\rfloor\ell}^{i\ell}\right\}.\]
    Namely, the net $N_{i\ell}(x_0,x)$ is built from the ends of the wedge $\Ss_{i\ell}(x_0,x)$ to the pointer (the pointer is not necessarily included in $N_{i\ell}(x_0,x)$), see Figure \ref{fig:net}. 
    Since each convex path $P_1^{i\ell}$ and $P_2^{i\ell}$ has length at most $i\ell$, in all cases $|N_{i\ell}(x_0,x)|\leq 2i+1$. 
    \begin{figure}
        \centering
        \includegraphics[width=0.25\linewidth]{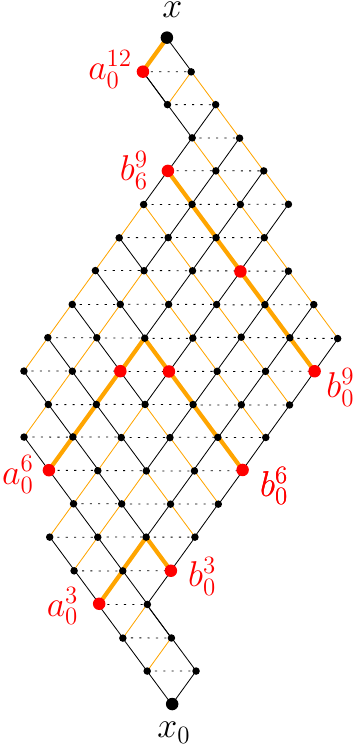}
        \caption{A net when $\ell=3$.}
        \label{fig:net}
    \end{figure}
    Finally, we define
    \[C_x = \bigcup_{i=0}^m N_{i\ell}(x_0,x).\]

    For every $z\in I(x_0,x)\setminus \{x_0\}$, let $0\leq i_z\leq m$ be the largest index such that $i_z\ell < d_{nor}(x_0,z)$, and let $i_{x_0}$=0.
    Then we define $p_x(z)=B_{3\ell}(z) \cap N_{i_z\ell}(x_0,x) \cap I(x_0,z)$.
    Since $B_{3\ell}(z)$ and $I(x_0,z)$ are convex, their intersection with the convex paths $P^{i_z\ell}_1(x_0,x)$ and $P^{i_z\ell}_2(x_0,x)$ consists of two convex subpaths $P'_1$ and $P'_2$ of length at most $6\ell$ each.  
    Hence $p_x(z)$ is the intersection of $N_{i_z\ell}(x_0,x)$ with $P'_1$ and $P'_2$.
    By maximality of $i_z$, the distance between $z$ and $\Ss_{i_z\ell}(x_0,x)$ is at most $2\ell$, thus the length of at least one of the paths $P'_1$ and $P'_2$ is at least $\ell$, or $P'_j = P^{i_z\ell}_j(x_0,x)$ for at least one $j\in\{1,2\}$.
    From the construction of the $\ell$-net $N_{i_z\ell}(x_0,x)$ it follows that $p_x(z)$ is nonempty. 

    Now we prove that the sets $C_x$ and the maps $p_x$ satisfy the conditions (a), (b), and (c). 
    Condition (b) directly follows from the definition of $p_x(z)$ as $B_{3\ell}(z) \cap N_{i_z\ell}(x_0,x) \cap I(x_0,z) \subseteq B_{3\ell}(z)$. 
    For condition (a), note that for $z\in I(x_0,y)$ we have $N_{i_z\ell}(x_0,y) = N_{i_z\ell}(x_0,x) \cap I(x_0,y)$.
    So, $C_x \cap I(x_0,y)=C_y$.
    Also, since   $I(x_0,z)\subseteq I(x_0,y)$,
    \[p_y(z)=B_{3\ell}(z) \cap N_{i_z\ell}(x_0,y) \cap I(x_0,z) = B_{3\ell}(z) \cap N_{i_z\ell}(x_0,x) \cap I(x_0,y) \cap I(x_0,z) = p_x(z).\]
    Hence, $p_x\vert_{I(x_0,y)}=p_y$.

    To show that (c) holds, first note that $B_{6\ell}(z)$ intersects at most $13$ of  the spheres $\Ss_{0}(x_0,x),\Ss_{\ell}(x_0,x),\Ss_{2\ell}(x_0,x),\ldots,\Ss_{m\ell}(x_0,x)$. 
    Each of the two convex paths of any of those spheres intersects $B_{6\ell}(z)$ in a subpath of length at most $12\ell$. Each such subpath contains at most $13$ vertices of the $\ell$-net of its sphere. So, the intersection of $B_{6\ell}(z)$ with each of the $13$ nets has size at most $26$.
    Hence, $|B_{6\ell}(z)\cap C_x|\leq 13\cdot 26=338$.
\end{proof}

Normal clique paths in systolic complexes have been defined in \cite{JaSw} to prove that systolic groups are biautomatic.
They are analogues of normal cube-paths in CAT(0) cube complexes.
For two vertices $u$ and  $v$ of a systolic complex $X$, a \emph{clique path} from $u$ to $v$ is a sequence of cliques $(\{ u\} = \sigma^{(0)}, \sigma^{(1)}, \ldots, \sigma^{(m)}=\{v\})$ such that any two consecutive cliques $\sigma^{(i)}$ and $\sigma^{(i+1)}$ are disjoint and their union induces a clique.
The \emph{normal clique path} $\gamma_{(u,v)}$ from $u$ to $v$ with $k=d(u,v)$ is a clique path defined as follows: for each $i$ running from $k$ to 0, we inductively define  $\sigma_{(u,v)}^{(i)}$ by setting $\sigma_{(u,v)}^{(k)}:=\{ v\}$ and $\sigma_{(u,v)}^{(i)}:=B^*_1(\sigma_{(u,v)}^{(i+1)}))\cap B_{i}(u)$. 
Notice that the normal clique paths are directed, since the cliques of $\gamma_{(u,v)}$ and $\gamma_{(v,u)}$ are not the same in general.
However both  $\gamma_{(u,v)}$ and $\gamma_{(v,u)}$ are subsets of the interval $I(u,v)$. 
In particular, if $I(u,v)$ is flat, then all cliques of $\gamma_{(u,v)}$ and $\gamma_{(v,u)}$ are vertices or horizontal edges of $I(u,v)$. 

Pick a base vertex $x_0$ and consider the normal clique paths $\gamma_{(x_0,y)}$, for a vertex $y$.
The following property of normal clique paths shows that if $x$ is close to $y$, then  $\gamma_{(x_0,y)}$ rapidly converges to the interval $I(x_0,x)$:

\begin{proposition} \cite{BlChGuOs} \label{p:convergence_of_normal_clique_paths} Let $X$ be a systolic complex with a base vertex $x_0$. If $x,y$ are two vertices of $X$ with $d(x,y)=k$ and $\gamma_{(x_0,y)}=(\{x_0\} =\sigma_0,\sigma_1,\ldots, \sigma_{m-1},\sigma_m=\{ y\})$
is the normal clique path from $x_0$ to $y$, then there exists an index $i\leq 2k$ such that the cliques $\sigma_0,\sigma_1,\ldots, \sigma_{m-i}$ are contained in the interval $I(x_0,x)$.
\end{proposition}

Now we will use Corollary~\ref{c:charANasdimfin}, Lemma~\ref{l:nets-2d}, and Proposition \ref{p:convergence_of_normal_clique_paths} to prove \Cref{thm:flat}. 

\begin{proof}[Proof of \Cref{thm:flat}]
    
For every $x\in X$, $\ell\in\mathbb{N}$, and $k\in\{0,\ldots,\ell\}$ we first define the auxiliary sets $\bar{S}(x,k,\ell)$.
For every $y\in B_k(x)$ consider the normal clique path $\gamma_{(x_0,y)}$ from $x_0$ to $y$ and let $\sigma_{(x_0,y)}^{(d(x_0,y)-2\ell)}$ be the clique at distance $2\ell$ from $y$ in $\gamma_{(x_0,y)}$ (or $\{ x_0\}$ if $d(x_0,y)<2\ell$).
We set \[\bar{S}(x,k,\ell)=\bigcup_{y\in B_k(x)}\sigma_{(x_0,y)}^{(d(x_0,y)-2\ell)}.\]
By Proposition \ref{p:convergence_of_normal_clique_paths},  $\bar{S}(x,k,\ell) \subseteq I(x_0,x)$. 
We define
\[S(x,k,\ell)=\bigcup p_x(\bar{S}(x,k,\ell))=\bigcup_{z\in \bar{S}(x,k,\ell)} p_x(z).\]
Finally, we show that the sets $S(x,k,\ell)$  verify the conditions (i)-(iv) in Corollary~\ref{c:charANasdimfin} with $D=6$. 
\begin{enumerate}[(i)]
  \item 
  By the definition of $\bar{S}(x,k,\ell)$, for any $z \in \bar{S}(x,k,\ell)$, we have $d(z,x)\leq 2\ell+k$.
  For such $z$ and any $w \in p_x(z)$ we know  by Lemma \ref{l:nets-2d}(b) that $d(z,w) \leq 3\ell$, which implies:
      \[S(x,k,\ell) \subseteq B_{5\ell+k}(x)\subseteq B_{6\ell}(x).\]
      
  \item For all $\ell \in \mathbb{N}$, $k,k'$ with $0 \leq k \leq k' \leq \ell$, and $x \in X$, we have $\bar{S}(x,k,\ell) \subseteq \bar{S}(x,k',\ell)$.
  Now immediately by the definition, $S(x,k,\ell) \subseteq S(x,k',\ell)$.
  
  \item Let $x'\sim x$ and $k=0,1,\ldots, \ell-1$. Since $B_k(x')\subseteq B_{k+1}(x)$, we 
  have $\bar{S}(x',k,\ell) \subseteq \bar{S}(x,k+1,\ell)$. To prove the inclusion $S(x',k,\ell) \subseteq S(x,k+1,\ell),$ we consider three cases. 

  \begin{case} $d(x_0,x) = d(x_0,x')+1$.
  \end{case}
  Then $x'\in I(x_0,x)$. By Lemma \ref{l:nets-2d}(a) $p_x|_{I(x_0,x')}=p_{x'}$, which implies
  \begin{align*}
  S(x',k,\ell)&=\bigcup p_{x'}(\bar{S}(x',k,\ell)) = \bigcup p_x(\bar{S}(x',k,\ell)) \\
  &\subseteq \bigcup p_x(\bar{S}(x,k+1,\ell))=S(x,k+1,\ell).
  \end{align*}

  \begin{case} $d(x_0,x) = d(x_0,x')-1$.
  \end{case}
  
   Then $x\in I(x_0,x')$. Since $\bar{S}(x',k,\ell)\subseteq \bar{S}(x,k+1,\ell)\subseteq I(x_0,x)$,  by Lemma \ref{l:nets-2d}(a), $p_{x'}|_{I(x_0,x)}=p_x$.  All this implies that 
   \begin{align*}
  S(x',k,\ell)&=\bigcup p_{x'}(\bar{S}(x',k,\ell)) = \bigcup p_x(\bar{S}(x',k,\ell)) \\
  &\subseteq \bigcup p_x(\bar{S}(x,k+1,\ell))=S(x,k+1,\ell).
  \end{align*}
   
  \begin{case} $d(x_0,x) = d(x_0,x')$.
  \end{case}
  
  Since $x\sim x'$, by the triangle condition, there exists $u \sim x,x'$ with $u\in I(x_0,x)\cap I(x_0,x')$. Since $\bar{S}(x',k,\ell)\subseteq \bar{S}(u,k+1,\ell)\subset I(x_0,u)$, we have
  \begin{align*}
  S(x',k,\ell)&=\bigcup p_{x'}(\bar{S}(x',k,\ell)) = \bigcup p_u(\bar{S}(x',k,\ell)) = \bigcup p_x(\bar{S}(x',k,\ell))\\
  &\subseteq \bigcup p_x(\bar{S}(x,k+1,\ell))=S(x,k+1,\ell),
  \end{align*}
  where the second equality follows from Lemma \ref{l:nets-2d}(a) and the inclusion follows from the inclusion $\bar{S}(x',k,\ell)\subseteq \bar{S}(x,k+1,\ell)$. This finishes the proof of  (iii). 

  \item By (i), $S(x,\ell,\ell) \subseteq B_{6\ell}(x)$.
  Hence, by definition, $S(x,\ell,\ell) \subseteq B_{6\ell}(x) \cap C_x$.
  Now by Lemma \ref{l:nets-2d}(c), we have $|S(x,\ell,\ell)| \leq 338$.
\end{enumerate}
\end{proof}

\section{Almost flat case}
\label{s:almflatproof}
In this section, we prove the \hyperref[thm:main]{Main Theorem} by generalizing the argument from Section~\ref{sec:flat}. Whenever possible, we refer to the proofs given there.
Let $u$ and $v$ be two vertices of a systolic complex $X$ with almost flat intervals with $d(u,v) = n$.
By Lemma \ref{l:structure-almost-flat-intervals}, there exist $0 = k_0< \ldots <k_m = n$ dividing the interval $I(u,v)$ into slices such that:
 \begin{itemize}
     \item[(1)]  the first  and the last slices  are chopped lozenges;
     \item[(2)] when moving from $u$ to $v$, chopped lozenges and zips alternate;
     \item[(3)] a chopped lozenge and a zip are either disjoint or are merged together along the common level, which is a single horizontal edge. 
 \end{itemize}

We now define the $k$-wedges $\Ss_k(u,v)$ of $I(u,v)$ as follows.
Add a ``ghost'' vertex $v_1$ to the top of $L_{k_0,k_1}(u,v)$, by joining it to both vertices in $L_{k_1}(u,v)$.
Then $L_{k_0,k_1}(u,v) \cup \{v_1\}$ is a burned lozenge, and we can define its wedges $\Ss_{m_0}(u,v),\ldots, \Ss_{m_1}(u,v)$ as in Section \ref{sec:flat}.
Now in $\Ss_{m_1}(u,v)$ we replace $v_1$ by both vertices in $L_{k_1+1}(u,v)$.

Next, we get to the first zip $L_{k_1,k_2}(u,v)$.
If $k_2-k_1 \leq 1$ we set $m_2=m_1$ and do nothing else.
Otherwise we define 
\[\Ss_{m_1+1}(u,v) = L_{k_1+2}(u,v), \ldots,\Ss_{m_2}(u,v) = L_{k_2}(u,v).\]

We now arrive to the second chopped lozenge $L_{k_2,k_3}(u,v)$.
Add again a ``ghost'' vertex $v_2$ to the top of $L_{k_2,k_3}(u,v)$.
Let $L_{k_2}(u,v) = \{u_2,u'_2\}$.
Inside of $L_{k_2,k_3}(u,v) \cup \{v_2\}$, consider the two burned lozenges $I(u_2, v_2)$ and $I(u'_2,v_2)$.
In each of them define wedges
\[\tilde\Ss_{0}(u_2,v_2) = \{u_2\},\ldots, \tilde\Ss_{\tilde{m}_3}(u_2,v_2) \text{ and } \tilde\Ss'_{0}(u'_2,v_2)= \{u'_2\},\ldots, \tilde\Ss'_{\tilde{m}'_3}(u'_2,v_2)\]
respectively.
Set $m_3 = m_2+\min\{\tilde{m}_3,\tilde{m}'_3\}$, and for $1 \leq j \leq \min\{\tilde{m}_3,\tilde{m}'_3\}$, set
\[\Ss_{m_2+j}(u,v) = ( \tilde\Ss_{j}(u_2,v_2) \setminus \cup_{s=0}^{j-1} \tilde\Ss'_s(u'_2,v_2)) \cup ( \tilde\Ss'_j(u'_2,v_2) \setminus \cup_{s=0}^{j-1} \tilde\Ss_s (u_2,v_2) ).\]
Finally, in $\Ss_{m_3}(u,v)$, replace $v_2$ by both vertices in $L_{k_3+1}(u,v)$.

\begin{figure}
    \centering
    \includegraphics[width=0.6\linewidth]{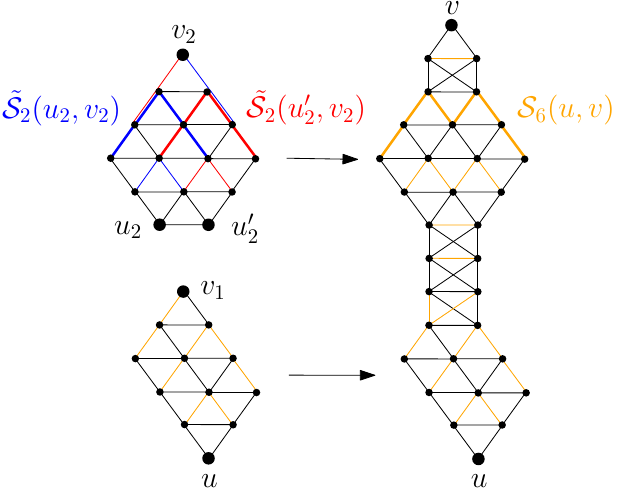}
    \caption{Wedges in an almost flat interval.}
    \label{fig:wedge-almostflat}
\end{figure}

Continuing with this process we define all wedges $\Ss_k(u,v)$ of $I(u,v)$ (note that in the last chopped lozenge, we do not add a ghost vertex and just use $v$). See Figure \ref{fig:wedge-almostflat} for an example. 

As in the proof of Theorem~\ref{thm:flat} the central point in the current proof is the following analogue of Lemma~\ref{l:nets-2d}.

\begin{lemma}\label{l:nets-almostflat}
Let $X$ be a systolic complex with almost flat intervals and let $x_0\in X$ be a vertex. Then for all $\ell \in \mathbb{N}$ and $x \in X$, there are subsets $C_x \subseteq I(x_0,x)$, and  maps $p_x\colon I(x_0,x) \rightarrow \mathcal{P}(C_x)\setminus \{ \emptyset\}$, where $\mathcal{P}(C_x)$ denotes the power set of $C_x$, satisfying:
\begin{enumerate}[(a)]
  \item If $d(x,y)=1$ and $y \in I(x_0,x)$, then $C_x \cap I(x_0,y)=C_y$, and $p_x\vert_{I(x_0,y)}=p_y$;
  \item For all $z \in I(x_0,x)$ and $w \in p_x(z)$, we have $d(z,w) \leq 3\ell$; 
  \item For any $z\in I(x_0,x)$, we have $|B_{6\ell}(z) \cap C_x| \leq 676$.
\end{enumerate}
\end{lemma}

\begin{proof}
    We proceed as in the proof of Lemma \ref{l:nets-2d}, by first defining the sets $C_x$ for every $x\in X$ and $\ell \in \mathbb{N}$.
    Assume $I(x_0,x)$ is partitioned into wedges $\Ss_0(x_0,x),\ldots,\Ss_M(x_0,x)$.
    Consider the wedges $\Ss_{0}(x_0,x),\Ss_{\ell}(x_0,x),\Ss_{2\ell}(x_0,x),\ldots,\Ss_{m\ell}(x_0,x)$ where $m=\lfloor\frac{M}{\ell}\rfloor$.

    We define an $\ell$-net $N_{i\ell}(x_0,x)$ for each  wedge $\Ss_{i\ell}(x_0,x)$.
    There are three cases to consider depending on whether the respective wedge is in the first chopped lozenge, in a zip, or in another chopped lozenge.

    \begin{case} $\Ss_{i\ell}(x_0,x)$ is in the first chopped lozenge.
    \end{case}
    Then we consider the burned lozenge obtained by adding the ghost vertex $v_1$ as before.
    Then we define $N_{i\ell}(x_0,x)$ just as in the proof of Lemma \ref{l:nets-2d}.
    In case $v_1$ belongs to the net, we replace it by both vertices in $L_{k_1+1}$.

    \begin{case} $\Ss_{i\ell}(x_0,x)$ is in a zip.
    \end{case}
    Then the wedge consists of just two vertices and we set $N_{i\ell}(x_0,x) = \Ss_{i\ell}(x_0,x)$.

    \begin{case} $\Ss_{i\ell}(x_0,x)$ is in another chopped lozenge.
    \end{case}
    Suppose that $\Ss_{i\ell}(x_0,x)$ is in $L_{k_h,k_{h+1}}$ and let $L_{k_h} = \{u,u'\}$.
    As before we add a ghost vertex $v$ and obtain two burned lozenges $I(u,v)$ and $I(u',v)$.
    Then the wedge is of the form
    \[\Ss_{i\ell}(x_0,x) = \big( \tilde\Ss_j(u,v) \setminus \cup_{s=0}^{j-1} \tilde\Ss'_s (u',v) \big) \cup \big( \tilde\Ss'_j(u',v) \setminus \cup_{s=0}^{j-1} \tilde\Ss_s (u,v)\big).\]
    We set
    \[N_{i\ell}(x_0,x) = \big( \tilde{N}_j(u,v) \setminus \cup_{s=0}^{j-1} \tilde\Ss'_s (u',v) \big) \cup \big( \tilde{N}_j(u',v) \setminus \cup_{s=0}^{j-1} \tilde\Ss_s (u,v)\big),\]
    where the $\tilde{N}_j(u,v)$ and $\tilde{N}_j(u',v)$ are the corresponding nets in $I(u,v)$ and $I(u',v)$ that are constructed as in the flat case. 

    Finally, we define
    \[C_x = \bigcup_{i=0}^m N_{i\ell}(x_0,x).\] See Figure \ref{fig:net-almostflat} for an illustration.

    \begin{figure}
        \centering
        \includegraphics[width=0.65\linewidth]{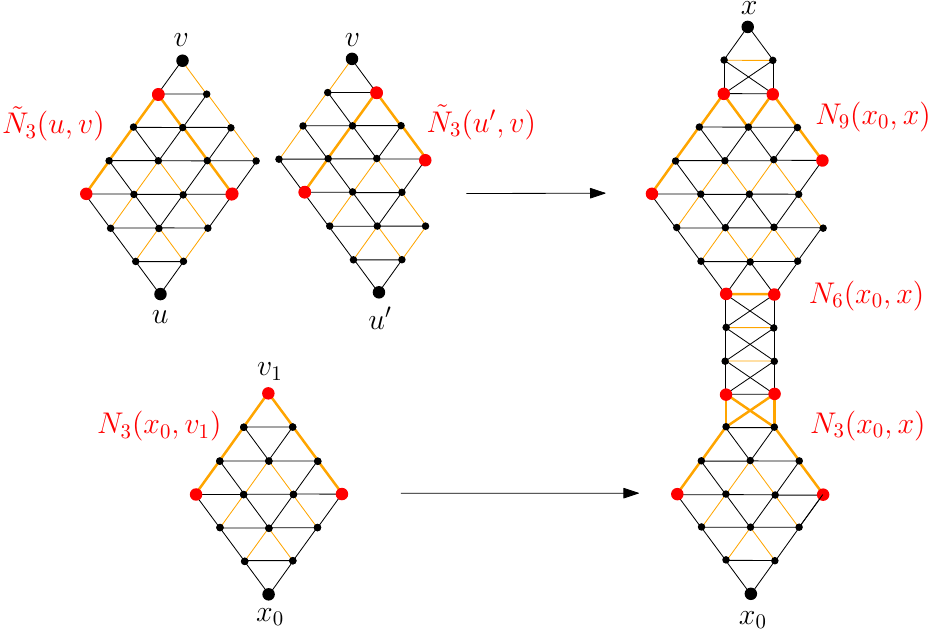}
        \caption{Right: an example of $C_x$ with $\ell = 3$ in an almost flat interval. Left to right: illustration of how nets in chopped lozenges are obtained from nets in corresponding flat intervals.}
        \label{fig:net-almostflat}
        
    \end{figure}

    Let $z\in I(x_0,x)$.
    We know $z$ belongs to some wedge $\Ss_p(x_0,x)$.
    Let $i_{x_0}=0$, and for $z\neq x_0$, let $0\leq i_z\leq m$ be the largest index such that $i_z\ell < p$.
    Then we set $p_x(z)=B_{3\ell}(z) \cap N_{i_z\ell}(x_0,x) \cap I(x_0,z)$.
    Note that $d(z,\Ss_{i_z\ell}(x_0,x))\leq 2\ell$, so as in the proof of Lemma \ref{l:nets-2d}, $p_x(z)$ is nonempty.

    Now we prove that the sets $C_x$ and the maps $p_x$ satisfy conditions (a), (b), and (c). 
    Condition (b) directly follows from the definition of $p_x(z)$ as $B_{3\ell}(z) \cap N_{i_z\ell}(x_0,x) \cap I(x_0,z) \subseteq B_{3\ell}(z)$. 
    For condition (a), note that for $z\in I(x_0,y)$ we have $N_{i_z\ell}(x_0,y) = N_{i_z\ell}(x_0,x) \cap I(x_0,y)$.
    So, $C_x \cap I(x_0,y)=C_y$.
    Also, since   $I(x_0,z)\subseteq I(x_0,y)$, we get 
    \[p_y(z)=B_{3\ell}(z) \cap N_{i_z\ell}(x_0,y) \cap I(x_0,z) = B_{3\ell}(z) \cap N_{i_z\ell}(x_0,x) \cap I(x_0,y) \cap I(x_0,z) = p_x(z).\]
    Hence, $p_x\vert_{I(x_0,y)}=p_y$.

    To show that (c) holds, first note that $B_{6\ell}(z)$ intersects at most $12$ of  the wedges $\Ss_{0}(x_0,x),\Ss_{\ell}(x_0,x),\Ss_{2\ell}(x_0,x),\ldots,\Ss_{m\ell}(x_0,x)$.
    In the worst case scenario, each of those wedges is the union of two flat-case wedges.
    So by the same argument as in the proof of Lemma \ref{l:nets-2d} we now get $|B_{6\ell}(z)\cap C_x|\leq 2\cdot 338 = 676$.
\end{proof}

\begin{proof}[Proof of {\hyperref[thm:main]{Main Theorem}}] We proceed exactly in the same way as in the proof of Theorem~\ref{thm:flat}, replacing Lemma \ref{l:nets-2d} by Lemma \ref{l:nets-almostflat}.
\end{proof}

\begin{proof}[Proof of Corollary \ref{mcor:gc3t6}] By tameness, the 
graphical $C(3)$--$T(6)$ complex $X_{(\Omega,\Theta)}$ (with the $1$--skeleton metric) is quasi-isometric to $\Th(X_{(\Omega,\Theta)})$, which is a systolic complex with almost flat intervals by \cite[Propositions 2.12 \& 2.13] {BlChGuOs} (see Section~\ref{s:c3t6} for details). Now, the result follows from {\hyperref[thm:main]{Main Theorem}}.
\end{proof}

\begin{remark}
\label{rem:nottame}
There exist locally finite non-tame simply connected graphical $C(3)$--$T(6)$ small cancellation complexes of infinite asymptotic dimension (with respect to the path metric on their $1$--skeleta).
For example, consider a family $\Omega=(\Omega_i)_{i=1}^{\infty}$ where all $\Omega_i$ are finite regular graphs of degree $3$, with $\mr{girth}(\Omega_i)\to\infty$, as $i\to \infty$. 
Let $X_{(\Omega,\Theta)}$ be the a graphical complex given by the natural embedding $f\colon \Omega \to \Theta$, where $\Theta$ is a connected graph obtained from the disjoint union of the $\Omega_i$ by adding an edge connecting a vertex of $\Omega_j$ with a vertex of $\Omega_{j+1}$, for every $j$.
By \cite[Theorem 1.2]{Will2011}  the asymptotic dimension of $\Theta$, and hence of $X_{(\Omega,\Theta)}$ (with the one skeleton path metric) is infinite. A CAT(0) version of this example appears in \cite{Kon2012, FukMat2026}.
\end{remark}

\begin{proof}[Proof of Corollary~\ref{mcor:c3t6group}]
    Since pieces in graphical $C(3)$--$T(6)$ complexes have length $1$, graphical $C(3)$--$T(6)$ presentations with finitely many generators (that is, the set $S$ as in Section~\ref{s:c3t6} is finite) have to have finitely many relations (that is, $\Omega$ is finite). Hence, there is a universal bound on the size of `relators' (that is, finite graphs $\Omega_i$). Therefore, the Cayley graph of the group is quasi-isometric to its thickening, and the result follows from {\hyperref[thm:main]{Main Theorem}}.
\end{proof}

\begin{remark}\label{rem:infgengr}
    Without the assumption of finite generation we can construct a counterexample following the one in Remark~\ref{rem:nottame}. Let $\Omega$ be a family from there. Let $S$ be a countable infinite set and consider an $S\cup S^{-1}$--labelling $f$ of edges of $\Omega$ such that no two edges are labelled by the same element or by inverse elements. Then the presentation 
    $\mathcal{P}=\langle S:f\rangle$ is an infinite graphical $C(3)$--$T(6)$ presentation and the corresponding Cayley graph contains isometrically embedded copies of $\Omega_i$ hence has infinite asymptotic dimension, by \cite[Theorem 1.2]{Will2011}.
\end{remark}

\section{Final remarks and open questions}
\label{s:final}

The question of finiteness of the asymptotic dimension is open for CAT(0) groups. Even the $2$-dimensional case---the most relevant for the current paper---is not known:

\begin{question}
	Have groups acting geometrically on $2$-dimensional CAT(0) complexes finite asymptotic (Assouad--Nagata) dimension?
\end{question}

Note that simply connected graphical $C(3)$--$T(6)$ small cancellation complexes can be equipped with a CAT(0) metric by \cite{Gui2025}.

\begin{question}
	Have classical $C(6)$ and $C(4)$--$T(4)$ small cancellation groups finite asymptotic (Assouad--Nagata) dimension?
\end{question}

Note that there are finitely generated but infinitely presented graphical $C(6)$ small cancellation groups with infinite asymptotic dimension (even without Yu's property A, and not coarsely embeddable into a Hilbert space) \cite{Osa-scl}.

\begin{question}
	Have finitely presented graphical $C(6)$ and $C(4)$--$T(4)$ small cancellation groups finite asymptotic (Assouad--Nagata) dimension?
\end{question}

Of course, we believe that the bound of the Assouad--Nagata dimension we got in the \hyperref[thm:main]{Main Theorem} is far from being optimal. We believe that the answer to the following question is affirmative.

\begin{question}
	Have tame graphical $C(3)$--$T(6)$ small cancellation complexes and $2$-dimensional systolic complexes asymptotic (Assouad--Nagata) dimension $2$?
\end{question}

Finally, this article gives a partial answer to the following question(s), to which we believe the answer is positive in general (cf.\ \cite[Remark 4.12]{EnMa}):

\begin{question}
	Have (fit) systolic groups finite asymptotic (Assouad--Nagata) dimension?
\end{question}

\appendix
\section{Buildings have finite asymptotic Assouad--Nagata dimension}
\label{s:app}

In \cite{DymSchi2009} Dymara and Schick prove that the asymptotic dimension of a building associated to a Coxeter group is equal to the asymptotic dimension of the group, and hence finite. We observe that their proof gives also the analogous bound for the asymptotic Assouad--Nagata dimension. 
We follow the notation and the layout of the proof from  \cite{DymSchi2009}, copying large parts from there and notifying only the important modifications, concerning quantitative estimates of the constants relevant in the new setting.
\medskip

Here $(W,S)$ is a finitely generated Coxeter system, $\mathfrak C$ is a building 
with Weyl group $W$, $|{\mathfrak C}|$ is the Davis realization of $\mathfrak C$, and $|W|$ denotes the Davis complex.
We fix a chamber $B\in {\mathfrak C}$ and define the $B$-based folding map
as $\pi\colon {\mathfrak C}\to W$ (also $\pi\colon |{\mathfrak C}|\to |W|$), with $\pi(B)=1$.
We consider two metric spaces: the discrete metric space $({\mathfrak C},\delta)$ with the gallery metric, and the $CAT(0)$
metric space $(|{\mathfrak C}|,\mathrm{d})$. These spaces are quasi-isometric, so there exists a constant $\cona>1$, depending only on $|W|$, such that
for every two chambers $c, c'$, and any two points $x\in |c|, x'\in |c'|$ we have
\begin{align}
	\label{f:QI}
	\cona^{-1} \mathrm{d}(x,x')-\cona \leq \delta (c,c')\leq \cona \mathrm{d}(x,x')+\cona.
\end{align}
For any set $U\subset {\mathfrak C}$ we will speak about its $\delta$- and $\mathrm{d}$-diameters of $U$ and denote them respectively by $\diam_\delta(U)$ and $\diam_{\mathrm{d}}(U)$. $W$ is isometrically embedded in $\mathfrak C$,  and
$|W|$ is isometrically embedded in $|\mathfrak C|$.  Thus we will denote the
distance functions of $W$ and $|W|$ also by $\delta$ and $\mathrm{d}$, respectively.

We use the following definition of the asymptotic Assouad--Nagata dimension:
A metric space $(X,d)$ has asymptotic Assouad--Nagata dimension $\le n$ if there exist constants $\cond>0$ and $D_0>1$ such that  
for any  $D>D_0$ there exist $n+1$ families 
$\U ^0,\ldots,{\U} ^n$ of subsets of $X$ such that: (1) $\bigcup_i{\U} ^i$ is a  
cover of $X$ consisting of sets of diameter at most $\cond D$, and (2) for every $i$ any two sets $U,U'\in{\U} ^i$ are $D$-disjoint:
$d(U,U'):=\inf\{d(x,y): x\in U, y\in U'\}\ge D$.
This definition is equivalent to the one we give in Definition \ref{d:as-Assouad-Nagata}.

Januszkiewicz \cite{Janu2002} showed that every Coxeter group admits a proper isometric action on a finite product of trees, hence the asymptotic Assouad--Nagata dimension of every Coxeter group is finite.

\begin{thmapp}
	\label{t:app}
	The asymptotic Assouad--Nagata dimension of any building associated with a Coxeter group $W$ is finite and equal to the asymptotic Assouad--Nagata dimension of $W$.
\end{thmapp}

\begin{proof}[Proof of Theorem~\ref{t:app}]
	
	The reader is advised to read the following while reading simultaneously the proof of \cite[Theorem 1]{DymSchi2009}. 
	
	Let $n$ be the asymptotic Assouad--Nagata dimension of the apartment $W$ of $\mathfrak C$, in particular we have the constants $\cond$ and $D_0$ as in the definition above. Because $W$
	embeds isometrically into $\mathfrak C$, the asymptotic Assouad--Nagata dimension of $\mathfrak C$ is at least
	$n$. We have to prove that its asymptotic Assouad--Nagata dimension is $\le n$. 
	
	Choose therefore $D>D_0$.
	Let $\U^0,\ldots,\U^n$ be families of $2D$-disjoint sets in $|W|$ of $\mathrm{d}$-diameter at most $2\cond D$ 
	such that $\bigcup_i\U^i$ is a cover of $|W|$. For  
	$U\in \U^i$ let $N_D(U)=\{x\in |W|: \mathrm{d}(x,U)<D\}$,
	and let $\conv{N_D(U)}$ be  the convex hull of this set in 
	the $CAT(0)$ space $|W|$. Note that 
	\begin{align}
		\label{f:0a}
		\diam_{\mathrm{d}}(\conv{N_D(U)})\le (2\cond+2)D.
	\end{align}
	Let ${\mathcal C}_U$ be the set of path--connected components of 
	$\pi^{-1}(\conv{N_D(U)})$. Put $\V_U=\{\pi^{-1}(U)\cap A: 
	A\in{\mathcal C}_U\}$ and $\V^i=\bigcup_{U\in\U^i}\V_U$.
	It is obvious that $\bigcup_i\V^i$ is a cover of $|\mathfrak C|$.
	
	It is shown in \cite{DymSchi2009} that $\V^i$ are  families
	of $D$-disjoint sets (the proof uses only the $2D$-disjointness of $\U^i$).
	We show that diameters of sets in $\V^i$ are bounded by $\cone D$, for a constant $\cone >0$ depending only on $|W|$.
	
	Let $\kappa$ be the $\mathrm{d}$-diameter of the realization of a chamber. For $A\in {\mathcal C}_U$ we have that 
	$T(A):=\bigcup\{Res(p): p\in A\}\subseteq\pi^{-1}\bigl(T(\conv{N_D(U)})\bigr)$, and any two points in $A$ belong to chambers which can be connected
	by a gallery in $T(A)$. Hence the set
	$A$ is contained in a gallery--connected component of $\pi^{-1}\bigl(T(\conv{N_D(U)})\bigr)$. We have that the $\mathrm{d}$-diameter of $T(\conv{N_D(U)})$ is, via (\ref{f:0a}), bounded by
	\begin{align}
		\label{f:00}
			\diam_{\mathrm{d}}\bigl(T(\conv{N_D(U)})\bigr)\leq (2\cond+2)D+2\kappa\leq (2\cond+2+2\kappa)D.	
	\end{align}
	
	The required bound on the elements of $\V^i$ follows from the following lemma, which is a slight (quantitative) modification of \cite[Lemma 1]{DymSchi2009}.
	\begin{lemma}
		\label{l:app1}
		There exists a constant $\conb$ depending only on $|W|$ such that for any $N\geq 1$ if $U$ is a subset of $W$ of 
		$\delta$-diameter $\le N$, then any gallery--connected component $V$
		of $\pi^{-1}(U)$ has $\delta$-diameter $\le \conb N$.
	\end{lemma}
	
	Assuming Lemma~\ref{l:app1} holds, we conclude the proof of Theorem~\ref{t:app} by setting the constant $\cone:=\cona \conb \cona (2\cond+2+2\kappa) + \cona \conb \cona  +\cona^2 $. It is so because, by (\ref{f:00}), for every gallery--connected component $V$
	of $\pi^{-1}\bigl(T(\conv{N_D(U)})\bigr)$ we have, using (\ref{f:QI}):
	\begin{align*}
		\diam_{\mathrm{d}}(V)\leq & \cona \diam_\delta(V) +\cona^2 
		\leq \cona \conb \diam_\delta \bigl(T(\conv{N_D(U)})\bigr) +\cona^2  \\
		\leq &\cona \conb \cona \diam_{\mathrm{d}} \bigl(T(\conv{N_D(U)})\bigr) + \cona \conb \cona +\cona^2 \\
		\leq &\cona \conb \cona (2\cond+2+2\kappa)D + \cona \conb \cona D +\cona^2 D.
	\end{align*}
	\end{proof}		
	The proof of Lemma~\ref{l:app1} follows the proof of \cite[Lemma 1]{DymSchi2009}, where we use the following slight (quantitative) modification of \cite[Lemma 2]{DymSchi2009}.
	
	\begin{lemma}
		\label{l:app2}
		Let $W$ be a Coxeter group. There exists a constant $\conc$ depending only on $|W|$ with the following property.
		For any $R>0$ and any subset $U\subset |W|$ of diameter
		$R$ satisfying  ${\mathrm d}(|1|,U)>\conc R + \conc$ there exists a codimension-one face of $|1|$ such that 
		the wall containing that face separates $|1|$ from $U$.
	\end{lemma}
	
	\begin{proof}[Proof of Lemma~\ref{l:app2}]
		The constant $\conc$ is given in the proof of \cite[Lemma 2]{DymSchi2009}, as follows: Pick $r$ greater than twice the diameter of a chamber. Let $\epsilon$ be the Lebesgue number of the finite open cover $\{(M_s^-\setminus B_r(b))\cap S_r(b): s\in S\}$ of the (compact) sphere  
		$S_r(b)$, where $M_s^-$ is the open half-space with boundary $M_s$ 
		that does not intersect $|1|$, for $M_s$ being the wall containing the $s$-face of the chamber $|1|$ and $b$ being the barycenter of $|1|$.
		Then we can take $\conc=\max \{r/\epsilon, r\}$
        \footnote{The argument in the proof of \cite[Lemma 2]{DymSchi2009} hinges on the fact that $U\subseteq  |W|\setminus B_r(b)$. Therefore one needs the additive term in ${\mathrm{d}}(|1|,U)>\conc R + \conc$ to deal with small $R>0$.}.
	\end{proof}
	
	\begin{proof}[Proof of Lemma~\ref{l:app1}]
		Let $U$ be a subset of $W$ of $\delta$-diameter $\le N$. Then the $\mathrm{d}$-diameter of $N(U)$ (the $1$-neighbourhood of $U$ in the $\delta$-metric) is, by (\ref{f:QI}),
		\begin{align}
			\label{f:01}
			\diam_\mathrm{d} |N(U)|\leq \cona N+ \cona^2 +2\kappa.
		\end{align}
		Iterated application of Lemma~\ref{l:app2} provides a minimal gallery $\gamma=(1,w_1,\ldots,w_k)$ 
		such that the wall between $w_i$ and $w_{i+1}$ separates $w_i$ from $N(U)$ 
		and 
		\begin{align}
			\label{f:02}
			\mathrm{d}(|w_k|,|N(U)|)\le \conc (\cona N+\cona^2 +2\kappa) + \conc.
		\end{align}
		Now, every chamber which
		belongs to $U$ can be joined to $1$ by a minimal gallery which 
		is a concatenation
		$\gamma\sigma$, i.e.\ extending $\gamma$. By \cite[Lemma 3]{DymSchi2009}, for every chamber $c$ of $V$  there is a minimal gallery from $B$ to $c$ (in $\mathfrak C$) whose $\pi$-projection
		extends $\gamma$, and there exists a chamber $e\in\pi^{-1}(w_k)$ such that any such minimal gallery passes
		through $e$. If $\alpha$ is the sub-gallery of such a gallery between $c$ and $e$ then the length of $\alpha$ is the same as the length of its image $\pi(\alpha)$ hence, by  (\ref{f:QI}), (\ref{f:01}), and (\ref{f:02}), we have
		\begin{align*}
			|\alpha|\leq& \cona (\mathrm{d}(|w_k|,|N(U)|) + \diam_\mathrm{d} (N(U))) +\cona \\
			 \leq &\cona \conc (\cona N+\cona^2 +2\kappa) + \cona\conc+ \cona (\cona N+ \cona^2 +2\kappa) +\cona. 
		\end{align*} 
		Therefore the $\delta$-diameter $\diam_{\delta}(V)$ of $V$ is at most
		\begin{align*}
			\label{f:04}
			2\bigl((\cona \conc +\cona)(\cona N+\cona^2 +2\kappa)+\cona \conc +\cona  \bigr) \leq 2 (\cona \conc +\cona)(\cona +\cona^2 +2\kappa +1) N. 
		\end{align*}
		Hence, the constant $\conb= 2(\cona \conc +\cona)(\cona +\cona^2 +2\kappa +1)$ is as required.
		\end{proof}

\begin{rem-app}
	Lang and Schlichenmaier \cite{LaSchli2005} show that the Assouad--Nagata dimension of affine buildings coincides with the Assouad--Nagata dimension of their appartments. Their proof does not extend to the case of general buildings, but their result covers also non-discrete affine buildings.
	
	Matsnev \cite{Matsnev2005} proved independently that affine buildings have finite asymptotic dimension.
\end{rem-app}

\subsection*{Acknowledgements} We thank Jan Dymara and Urs Lang for valuable discussions. The work presented here was partially supported by the Carlsberg Foundation, grant CF23-1226. V.\ C.\ was partially supported by the ANR project MIMETIQUE (ANR-25-CE48-4089-01). D.\ O.\ was partially supported by the (Polish) Narodowe Centrum Nauki (NCN), project no.\ UMO-2025/59/B/ST1/03217.

\bibliographystyle{amsalpha}
\bibliography{mybib}

\end{document}